\documentclass[12pt, oneside]{amsart}

\usepackage{amssymb}
 \usepackage[foot]{amsaddr}
 
\usepackage{mathtools}
\usepackage{xcolor}
\usepackage{url} 
\usepackage[breaklinks=true]{hyperref}
\usepackage{enumitem}

\usepackage{booktabs}
\usepackage{array}
\usepackage{adjustbox}
\usepackage{multirow}
\usepackage{makecell}

\usepackage{graphicx}
\usepackage{subfig}

\usepackage{float}

\usepackage{todonotes}

\newtheorem{thm}{Theorem}
\newtheorem{prop}[thm]{Proposition}

\theoremstyle{definition}

\theoremstyle{remark}
\newtheorem{remark}[thm]{Remark}

\newcommand{\defn}{\ensuremath{:=}}

\newcommand{\N}{\mathbb{N}}
\newcommand{\R}{\mathbb{R}}

\newcommand{\cP}{\mathcal{P}}

\newcommand{\cU}{\mathcal{U}}

\newcommand{\bfi}{{\mathbf{i}}}

\newcommand{\set}[1]{\left\{#1\right\}}

\newcommand{\norm}[2][2]{\left\lVert#2\right\rVert_{#1}}

\newcommand{\inner}[2]{\left\langle#1,#2\right\rangle}

\newcommand{\abs}[1]{\left \vert #1 \right \vert }

\let\epsilon\varepsilon

\let\phi\varphi

\newcommand{\first}[1]{\textbf{#1}}

\newcommand*{\doi}[1]{\href{https://doi.org/#1}{doi: #1}}

\newcommand{\AlgSymbol}{{A}}
\newcommand{\LSSymbol}{\AlgSymbol^{\mathrm{LS}}}
\newcommand{\SmoSymbol}{\AlgSymbol^{\mathrm{S}}}
\newcommand{\Alg}[2]{#1_{#2}}
\newcommand{\LS}[1]{\Alg{\LSSymbol}{#1}}
\newcommand{\wLS}[1]{\Alg{A^{\mathrm{wLS}}}{#1}}
\newcommand{\SM}[1]{\Alg{\SmoSymbol}{#1}}
\newcommand{\Sm}[1]{\SM{#1}}

\newcommand{\widthoffigures}{1.0}

\begin{document}
\title[Sparse grids 
vs.~random points]{
Sparse grids 
vs.~random points \\
for high-dimensional \\ 
polynomial approximation
}

\author{Jakob Eggl$^1$, Elias Mindlberger$^1$ and Mario Ullrich$^{1,2}$}


\address{$^1$Department of Quantum Information and Computation, 
Johannes Kepler University Linz, Austria, 
$^2$Institute of Analysis, Johannes Kepler University Linz, Austria
}
\email{research.jakob.eggl@gmail.com} 
\email{research@mindlberger.com}
\email{mario.ullrich@jku.at}


\date{\today}

\begin{abstract} 
We study polynomial approximation on a 
$d$-cube, 
where $d$ is large, 
and compare interpolation on sparse grids, 
aka~Smolyak's algorithm (SA), 
with a simple least squares method based on 
randomly generated 
points (LS) using standard benchmark functions. 

Our main motivation is the influential paper 
[Barthelmann, Novak, Ritter: 
High dimensional polynomial interpolation on sparse
grids, Adv.~Comput.~Math.~12, 2000]. 
We repeat and extend their theoretical analysis and numerical experiments 
for SA and compare to LS  
in dimensions up to 100. 
Our extensive experiments 
demonstrate that LS, even with only slight oversampling, consistently matches the accuracy of SA in low dimensions. 
In high dimensions, 
however, LS shows clear superiority.
\end{abstract}


\maketitle


\section{Introduction}

The approximation of high-dimensional functions 
based on a finite number of function evaluations 
is a fundamental task in computational science. 
Among the most widely studied approaches, 
for functions defined on the $d$-cube $[0,1]^d$, 
are interpolation on sparse grids (SA) and least-squares (LS) methods 
based on general point sets. 
The advantage of SA, see Section~\ref{sec:SA} for 
details, 
is that it is tailored for functions with a certain tensor product structure 
leading to nice performance and theoretical error guarantees
for corresponding problems. 

The advantage of LS, 
see Section~\ref{sec:LS}, 
is its flexibility (and simplicity) in the sense that, in contrast to SA, one does not need point sets of a specific form and cardinality.
The latter is important in high dimensions because 
known 
constructions of ``good'' point sets often 
require an exponential in $d$ number of points. 
To overcome this shortcoming, one may use random points which are easy to generate, also in large dimensions.

However, the use of \emph{randomness} generally does not allow for theoretical a priori error guarantees, but only for probabilistic bounds. 
In recent years, and based on advances in random matrix theory, 
it has been observed that 
LS based on iid random points 
is, with high probability, 
a ``good''
and to some extent universal method in a deterministic worst-case sense. 
See Section~\ref{sec:LS} for details. 
In other words, if we generate 
$n$ random points in~$[0,1]^d$ once and fix them (i.e., store them on the computer),  
the corresponding LS 
should
be ``good'' 
in many cases.

This paper presents a systematic numerical comparison of 
this method and SA 
using a set of randomly generated standard benchmark test functions. 
We basically follow~\cite{BNR00} and also discuss theoretical bounds for a model class of high-dimensional functions. 
We then provide extended numerical tests and compare with 
two different least squares methods: One with uniformly distributed points and equal weights, 
and one with points/weights according to the Chebyshev distribution, see Section~\ref{sec:LS}. 
Both LS methods will project on the same polynomial space as SA and use twice the number of points. 

The numerical results, which are described in more detail in Section~\ref{sec:experiments}, can be summarized as follows:
\begin{itemize}
\itemsep=1mm
\item In every scenario, one of the 
LS
methods is as good as SA. 
\item $2n$ random points seem to suffice for this performance, which is not yet supported by theoretical error bounds.
\item Often an unweighted LS with uniformly distributed points works, also for approximation in the (estimated) uniform norm.
\item In high dimensions, 
least squares often outperforms Smolyak's algorithm a lot.
\item As expected, LS is
more \emph{stable} when confronted with \emph{noise}.
\end{itemize}

\medskip

We also discuss the running times and difficulties with implementations in Section~\ref{sec:discussion}.
Our 
code 
and results are available at 
\begin{center}
\url{https://github.com/th3lias/ApxKit}. 
\end{center}

\newpage

\section{Sparse grids \& Smolyak's algorithm}
\label{sec:SA}

Smolyak's algorithm, i.e., interpolation on sparse grids, 
has been of high theoretical and practical 
interest in the last decades. 
We refer to~\cite{BG04,DTU18,JMMV14} and references
therein for detailed treatments and~\cite{Smo63} for the original source. 
We follow the description of \cite{BNR00}, and recommend reading this article as a preparation, see also~\cite{NR96,NR99}.

We aim to recover a function \(f: [0, 1]^d \to \R,\, d \in \N,\)
from $n$ samples 
$f(x_1),\dots,f(x_n)$, 
where we assume that the $x_i$ are known. 
In this section, we describe \emph{interpolatory methods}, 
i.e., algorithms that exactly reproduce the function values at the points $x_i$. \\
For \(d=1\), the general form of 
such an algorithm 
is given by
\begin{equation}\label{eq:1dInterpolant}
	\cU(f) \defn \sum_{j=1}^n 
        f(x_j)\,
        a_j
\end{equation}
where \(a_j \in C[0, 1]\) are chosen such that \(a_j(x_i) = \delta_{ij},\, i,j \in \N\). \\
In order to extend this to the multivariate case, \(d > 1\),
assume the existence of a sequence of algorithms 
\(\cU^i\) in the form
of \eqref{eq:1dInterpolant} for point sets \[
	X_i \defn \set{x_1^i, \dots, x_{m_i}^i}
\]
and $n=m_i$, 
where \(m: \N \to \N\) specifies the cardinality of each point set $X_i$, $i \in \N$. 
To simplify notation, we set $m_0=0$ and \(\cU^0 \defn 0\). 

We then form the tensor products of the univariate interpolations by 
\begin{equation}\label{eq:DdInterpolation}
	\cU^{\bfi}(f) \defn \left(\bigotimes_{k=1}^d \cU^{i_k}\right) (f) = \prod_{k=1}^{d} \sum_{j_k = 1}^{m_{i_k}} f\left( x_{j_1}^{i_1}, \dots, x_{j_d}^{i_d} \right) \left( a_{j_1}^{i_1} \otimes \cdots \otimes a_{j_d}^{i_d} \right)
\end{equation}
with \(\bfi \defn (i_1, \dots, i_d) \in \N_0^d\). 
This 
needs \(\prod_{k=1}^d m_{i_k}\) function evaluations. \\
\emph{Smolyak's algorithm} 
is obtained by 
using suitable linear combinations of these building blocks as follows. 
We define a 
\emph{resolution parameter} \(q \in \N\) with $q\ge d$ and only consider multiindices for which \(\norm[1]{\bfi} \leq q\). 
Using the \emph{differences} of the univariate interpolants defined by 
\(\Delta^i \defn \cU^{i} - \cU^{i-1}\), \(i \in \N\), 
we define Smolyak's algorithm by 
\begin{equation*}
    \Sm{q, d} \defn \sum_{\norm[1]{\bfi} \leq q} \bigotimes_{k=1}^d \Delta^{i_k},
\end{equation*}
where we again set \(\Delta^0 \defn 0\).
Equivalently, we may write 
\begin{equation*}
	\Sm{q,d} = \sum_{q-d+1 \leq \norm[1]{\bfi} \leq q} (-1)^{q-\norm[1]{\bfi}} \binom{d-1}{q - \norm[1]{\bfi}} \; 
    \cU^{\bfi},
\end{equation*}
see e.g.~\cite[Lemma 1]{WW95}. 
To compute the interpolant \(\Sm{q,d}(f)\), one only needs to sample the function at the \emph{sparse grid}
\begin{equation*}\label{eq:sparseGrid}
	H(q,d) \defn \bigcup_{q-d+1 \leq \norm[1]{\bfi} \le q} \bigtimes_{k=1}^d X_{i_k}. 
\end{equation*}

It remains to specify a sequence of univariate point sets $(X_i)$ and corresponding functions $a^i_j$, $j=1,\dots,m_i$, see \eqref{eq:1dInterpolant} and \eqref{eq:DdInterpolation}.
Here, as in~\cite{BNR00}, 
we suggest using rules based on polynomial interpolation at
the extrema of the Chebyshev polynomials, which is an especially popular choice.
These are defined by \(x_1^i = 1/2\) for \(m_i = 1\) and, 
given \(m_i > 1\), we set
\begin{equation}\label{eq:chebyExtrema}
	x_j^i \,:=\, - \frac12 \cos \frac{(j-1)\pi}{(m_i-1)} + 1/2 \qquad \text{for} \quad 
    j = 1, \dots, m_i.
\end{equation}
\begin{remark}
    Usually, this point set is given in the interval \([-1,1]\).
    For coherence with later results, we scale to \([0, 1]\) via \(x \mapsto (x+1)/2\).
\end{remark}

The functions $a^i_j$ are specified by the demand that the univariate methods $\cU^i$ 
reproduce all polynomials of degree less than $m_i$ exactly. 
That is, $\{a_j^i\}_{j=1}^{m_i}$ is the Lagrange basis of $X_i$. 

For the cardinality \(m_i\) we choose a doubling rule 
\begin{equation*}\label{eq:gridCard}
    m_1 = 1 
    \quad\text{and}\quad
    m_i = 2^{i-1}+1 \quad\text{for}\quad i > 1,
\end{equation*}
which ensures that the successive 
grids are nested, i.e., $X_i \subset X_{i+1}$, 
and thus 
$H(q, d) \subset H(q+1, d)$.
The main advantage of nested point sets
is that the function evaluations used for computing $\cU^i(f)$ can be reused for $\cU^{i+1}(f)$. 
Moreover, the corresponding Smolyak's algorithm is interpolatory on the sparse grid. 
In fact, the following has been proven in~\cite[Proposition~6]{BNR00}.

\begin{prop}
	Assume \(X_1 \subset X_2 \subset \dots\) and that \(\cU^i(f)\) are interpolatory, 
    that is, $\cU^i(f)(x) = f(x)$ for any $f \in C[0, 1]$, $x \in X_i$ and $i\in\N$. 
    Then, 
    \[
	   \Sm{q,d}(g)(x) = g(x)
	\]
    for any \(g \in C[0,1]^d\) and \(x \in H(q, d)\).
\end{prop}

\medskip

Another crucial property of Smolyak's algorithm 
is that certain exactness of the univariate 
building blocks is reflected also in the $d$-variate method. 
Let \(\cP(k,d)\) denote the space of polynomials 
from \(\R^d\) to \(\R\) of total degree at most \(k\), and 
recall that the univariate $\cU^i$ are designed 
to be exact on $\cP(m_i-1,1)$.
It turns out, see~\cite[Proposition~3 \& Theorem~4]{BNR00}, that Smolyak's algorithm is exact on some ``nonclassical'' space of polynomials. 

\begin{prop} \label{prop:space}
$\Sm{q,d}$ as described above 
    is exact on \[
	E(q, d) \defn \sum_{\norm[1]{\bfi} = q} 
\cP(m_{i_1}-1,1) \otimes \cdots \otimes \cP(m_{i_d}-1,1).
	\]
In particular, \(\Sm{d+k,d}\) is exact 
on $\cP(k,d)$.
\end{prop}

That is, $\Sm{d+k,d}$ is exact for all polynomials of total degree \(k\), but also for \emph{some} polynomials of higher degree. 
We will use $E(q,d)$ also as \emph{approximation space} for the least squares method, see Section~\ref{sec:LS}.

The number $N(q,d)$ of function evaluations 
used by $\Sm{q,d}$ satisfies some favorable bounds for large $d$. 
It has been observed in~\cite[Lemma~2]{NR99} that for the above construction 
(with $m_1=1$ and $m_2=3$) we have
\begin{equation}\label{eq:Nqd}
N(d+k,d) \,:=\, \# H(d+k,d) 
\,=\, \dim E(d+k,d) 
\,\approx\, \frac{2^k}{k!} d^k 
\end{equation}
for fixed $k$ and $d\to\infty$, 
where $a_d\approx b_d :\Leftrightarrow \lim_{d\to\infty}\frac{a_d}{b_d}=1$. 
Note that this compares well to $\dim \cP(k,d)\approx d^k/k!$ 
which is the minimal number of 
function values needed by any method that is exact on $\cP(k,d)$.

\medskip

Let us finally discuss a theoretical error bound, as in~\cite{BNR00}, that justifies the usage of $\Sm{q,d}$ in high dimensions.
For this, we define the model class 
\begin{equation*}\label{eq:Fsd}
F_d^s \defn \set{f: [0, 1]^d \to \R \colon D^\alpha f \text{ continuous if } 
\alpha_i \le s \text{ for all } i}
\end{equation*}
equipped with the norm 
\[
\|f\|_{F_d^s} \,:=\, \max\left\{\|D^\alpha f\|_\infty\colon \alpha\in\N_0^d,\, \alpha_i\le s \right\}, 
\]
where $\|f\|_\infty:=\sup_{x}|f(x)|$.

This class is tailored for multivariate functions that show some \emph{independence} of variables. 
In particular, for univariate functions $f_1,\dots,f_d\in F_1^s$, 
we have $f_1\otimes\dots\otimes f_d\in F_d^s$ 
with 
$\|f_1\otimes\dots\otimes f_d\|_{F_d^s} = \|f_1\|_{F_1^s}\cdots\|f_d\|_{F_1^s}$. 

\medskip

The proof of an error bound is now based on a suitable combination of the univariate properties. 
For this, note that the 
univariate (linear) operator 
\(\cU^i\) as defined above is exact on 
\(\cP(2^{i-1},1)\). 
Hence, we obtain the general bound 
\[
    \norm[\infty]{f - \cU^i(f)} \,\le\, 
    (1+\Lambda_{i})\, 
    \inf_{p \in \cP(2^{i-1},1)} 
    \norm[\infty]{f-p}, 
\]
where \(\Lambda_{i}\) is the Lebesgue constant
\begin{equation*}\label{eq:lebesgueConstant}
    \Lambda_{i} \defn 
    \sup_{\norm[\infty]{f}=1} \norm[\infty]{\cU^i(f)}.
\end{equation*}
Here, \(\cU^i\) is an interpolating operator 
as in \eqref{eq:1dInterpolant} 
for the 
Chebyshev points specified by \eqref{eq:chebyExtrema}, 
in which case it is known that \[
    \Lambda_{i} < \frac{2}{\pi} \, 
    \ln (m_i-1)
    + C
\]
for some $C > 0$, 
see~\cite[Satz 4]{EZ66}.
Together with the well-known Jackson bound 
$\inf_{p \in \cP(n,1)} \norm[\infty]{f-p}\lesssim n^{-s}\,\|f\|_{F_1^s}$, 
we obtain the error bound
\[
\norm[\infty]{f - \cU^i(f)} \,\le\, 
c_{1,s}\cdot i\cdot 2^{-is}\,\|f\|_{F_1^s}
\]
for all $f\in F_1^s$.

Based on this, the following bounds 
have been proven in \cite{BNR00}.

\begin{thm}\label{thm:SA}
For $q,d,s\in\N$, we have 
\[
\norm[\infty]{f - \Sm{q,d}(f)} \;\le\; c_{d,s}\, n^{-s} \left(\log 
n\right)^{(s+2)(d-1)+1}\,
\|f\|_{F_d^s} 
\]
for all $f\in F^s_{d}$, 
where 
$n:=N(q,d)$ is the number of points used by $\Sm{q,d}$ 
and 
\(c_{d, s}\) is a positive constant only depending on \(d\) and \(s\).
\end{thm}

\medskip

This bound shows an order of convergence that is, apart from the logarithmic term, 
independent of $d$. 
This makes it a promising model for high-dimensional approximation problems.

\begin{remark} 
We basically repeat Remark~9 of~\cite{BNR00}. \\
The bound from Theorem~\ref{thm:SA} is not optimal, and one may achieve 
\[
\norm[\infty]{f - A_n(f)} \;\lesssim\; n^{-s} \left(\log 
n\right)^{(s+1)(d-1)}\,
\|f\|_{F_d^s} 
\]
for a Smolyak algorithm $A_n$ using at most $n$ samples with the univariate formulas $\cU^i$ in the above construction chosen to yield the (optimal) order $2^{-is}$, see~\cite{Smo63, Tem86, WW95}. 
The optimal constructions that we are aware of, like splines, depend on the smoothness of the space. 
In contrast, the algorithm introduced here works for all $s$. Such methods are often called \emph{universal}.
\end{remark}

\bigskip

\section{Least Squares Approximation}
\label{sec:LS}

We now consider the conceptually simpler 
 \emph{weighted least squares estimator} 
\begin{equation*}\label{eq:wLS}
A_{X,V}^w(f) \,:=\, 
\underset{g\in V}{\rm argmin}\, \sum_{i=1}^n 
w_i\, \vert f(x_i)-g(x_i) \vert^2,
\end{equation*}
which is specified by a finite-dimensional space 
$V\subset C(D)$, 
for some set~$D$, 
a set of points $X:=\{x_1,\dots,x_n\}\subset D$ 
and weights $w=(w_1,\dots,w_n)$ with $w_i>0$.
We simply write $A_{X,V}$ for constant weights.

Assuming that the points $X$ are such that 
$g\in V$ and $g|_X=0$ imply $g=0$, 
which can only hold if $n\ge\dim(V)$, 
we can find $A_{X,V}^w(f)$ by computing the (in this case unique) least squares solution of the 
overdetermined system $Mz = y$ 
with $M=(\sqrt{w_i}\, b_j(x_i))_{i,j}$ and 
\(y = \left(\sqrt{w_i}\,f(x_j)\right)_{j=1}^n\), 
where $b_1,\dots,b_m\in V$ is an arbitrary basis of $V$.
To be precise, 
we have $A_{X,V}^w(f)=\sum_{j=1}^k z^*_{j}\, b_j$
for $z^*:=\inf_{z \in \R^k} \norm[2]{Mz - y}$.
Under the same condition, we also obtain 
$A_{X,V}^w(g)=g$ for all $g\in V$.
However, contrary to Smolyak's algorithm, the least squares method is, in general, not interpolatory unless $n=\dim(V)$. 

The method of least squares is clearly heavily studied and applied to various problems. 
We refer to~\cite{CCMNT15,Gr19,GNZ20,Jo92,LT84} and references therein. 
From an (approximation) theoretical point of view, the method was mostly analyzed for approximation in $L_2$-norm, 
and for special $V$ such as polynomials. 
The attempt to examine its approximation performance in more generality and for other (e.g.~uniform) norms seems to be quite new, see~\cite{BCLSV11,CL08,Reichel86}. 
Recent advances show that least squares methods (with suitable $V,X,w$) allow for optimal (in the sense of order of convergence) error bounds in different settings, and even improve (theoretically) over more tailored approaches, see e.g.~\cite{BPV19,DC22,DKU,KPUU25-general}.

Several of these advances are based on \emph{random sampling points},  
which is the motivation of this paper.

\medskip

Concerning error bounds, 
let us first fix some notation. 
For a measure $\mu$ on $D$, we denote by $L_p(\mu):=L_p(D,\mu)$, 
$1\le p<\infty$, 
the space of $p$-integrable functions with the usual norm 
$\|f\|_{L_p(\mu)}^p:=\int_D |f|^p d\mu$, 
and with the typical modification for $p=\infty$. 
We will mostly consider $D=[0,1]^d$ equipped with 
the Lebesgue measure.

The following general result has recently been proven in~\cite[Theorem~17]{KPUU25-uniform}.
It shows 
that, 
for every $D$, $\mu$ and $V$, 
there exist point sets of size close to $\dim(V)$ that are simultaneously ``good'' 
for approximation in all $L_p$-norms.
The proof is based on the close connection between LS and suitable discretization problems, see Remark~\ref{rem:disc}.

\begin{thm}\label{thm:LS-Lp}
Let \(\mu\) be a Borel probability measure on a compact set $D\subset\R^d$ and 
$V_n$ be an $n$-dimensional subspace of 
$C(D)$. \\
Then, there are points 
\(X=\{x_1, \dots, x_{4n}\} \subset D\) 
such that for all \(f \in C(D)\) and all \(1 \leq p \leq \infty\) we have 
\[
\norm[L_p(\mu)]{f - A_{X,V_n}(f)} \;\le\; 84\cdot n^{(1/2-1/p)_+}\cdot \inf_{g \in V_n} \norm[\infty]{f-g}, 
\]
where $a_+ \defn \max(0, a)$ for $a \in \R$. 
\end{thm}

The result in this generality cannot be improved (up to constants), 
see~\cite[Section~5.3]{KPUU25-uniform}. 
For $p=2$, Theorem~\ref{thm:LS-Lp} has been proven in~\cite{T21} for least squares methods with weights, and in~\cite{BPV19} the case $p=\infty$ for polynomials is discussed 
in more detail.

\medskip

Theorem~\ref{thm:LS-Lp} for $p=\infty$ 
has a factor $\sqrt{n}$ on the right hand side, 
and therefore cannot reproduce the upper bound from Theorem~\ref{thm:SA}. 
It seems to be unclear for which $p$ a suitable least squares method outperforms Smolyak's algorithm for the classes $F_d^k$, 
or other Sobolev spaces,  
and we refer to~\cite{KPUU25-general} for a discussion. 

Also in the case $p=2$, we do not obtain interesting bounds from 
Theorem~\ref{thm:LS-Lp} for examples relevant to this paper. 
The reason is that we compare the $L_2$-approximation power of least squares with the \emph{best approximation} on $V_n$ measured in sup-norm. 
The latter often has a worse dependence on $n$, see e.g.~\cite{DTU18}. 
Moreover, results like Theorem~\ref{thm:LS-Lp} do not lead to feasible constructions. 

These issues have been discussed in a series of recent papers, 
see~\cite{DKU,KU21,KU21b} for $p=2$, 
\cite{CDL13,CM17,DC22} for an analysis of \emph{randomized algorithms}, 
and 
\cite{BSU23,CD23} for approaches to make the above \emph{constructive}.
We only discuss one part on the story here, and refer to 
the more comprehensive works~\cite{GNZ20,KKLT,KPUU25-general,SU23} 
for more details, extensions and history. 

In fact, if we assume that the algorithm is applied only to 
$f\in F$ for some ``nice'' class $F\subset C(D)$, 
then a suitable least squares method 
has the optimal order of convergence,
see~\cite{DKU}. 
Moreover, 
we may use a point set of
independent identically distributed points 
if we allow some \emph{oversampling}, see~\cite{KU21,KU21b}. 

\begin{thm}\label{thm:LS-L2}
Let \(\mu\) be a Borel probability measure on a compact set $D\subset\R^d$ and 
$(V_n)_{n\in\N}$ be a sequence of $n$-dimensional subspaces of 
$C(D)$, i.e., $\dim(V_n)=n$ for all $n$. 
Moreover, for $F\subset C(D)$, we 
assume that 
\begin{equation*}
\sup_{f\in F}\,\inf_{g \in V_n} \norm[L_2(\mu)]{f-g} 
\;\le\; C\, n^{-\alpha} (\log n)^\beta 
\end{equation*}
for some $\alpha>1/2$ and $\beta,C>0$. \\
Then, there are constants $b,c\in\N$, depending only on $\alpha,\beta,C$, 
such that for all $n\in\N$, 
there are points 
\(X_n=\{x_1, \dots, x_{bn}\} \subset D\) 
and weights $w\in\R_+^{bn}$ 
with  
\[
\norm[L_2(\mu)]{f - A^w_{X_n,V_n}(f)} \;\le\; c\cdot n^{-\alpha} (\log n)^\beta
\]
for all \(f \in F\).

This result remains true, 
with probability $1-\frac{1}{b n^{2}}$, 
if one chooses $X_n$ 
to consist of $b\, n\log(n)$ independent random points 
with distribution $\rho_n\cdot d\mu$ 
and $w_k:=\rho_n(x_k)^{-1}$, 
where $\rho_n$ is some probability density that depends only 
on $V_n$ and $F$.
\end{thm}

We refer to~\cite{KPUU25-general} and~\cite[Thm.~3.7]{SU23} 
for this version of the theorem.

\medskip

Roughly speaking, this result shows that a suitable weighted least squares estimator can catch up even with the best-approximation in~$L_2$ on the same subspace (up to constants; or some $\log$) in a deterministic worst-case setting whenever $F$ is such that the latter has an error smaller than $n^{-1/2}$. 
The density $\rho_n$ of the last theorem, although given explicitly in~\cite{KU21,KU21b}, is not very handy and may be impossible to sample from in general. 
However, in some cases, the density reduces to a simple one. 

Let us discuss this with a (general) example 
of tensor product spaces. \\
First, consider 
$D_1=[0,1]$ 
equipped with the Lebesgue measure $\lambda$, and let $d\mu_1:=\nu\cdot d\lambda$ for some probability density $\nu\colon [0,1]\to\R_+$. 
Hence, $\langle f,g\rangle_{\nu}=\int_{0}^1 f(x) \, \overline{g(x)}\, \nu(x) dx$ is the inner product in $L_2([0,1],\mu_1)$.
Now let $b_0,b_1,\dots$ be an orthonormal basis (of polynomials) of $L_2([0,1],\mu_1)$, 
such as 
\begin{itemize}
    \item \emph{trigonometric polynomials}: $b_k^T(x):=e^{2\pi{\rm i}k x}$ for $\nu(x)=1$, or 
    \item \emph{Legendre polynomials}: $b_k^L$ for $\nu(x)=1$, or
    \item \emph{Chebyshev polynomials}: $b_0^C=1$ and 
    $$b^C_k(x) := \sqrt{2}\,\cos(k \arccos(2x-1)), \qquad k\ge1$$
    for
    $\nu(x)=\left(\pi \sqrt{x(1-x)}\right)^{-1}$ 
    being the \emph{Chebyshev weight}. 
\end{itemize}
Note that the polynomials are usually considered on $[-1,1]$.

Given such a sequence, we define  
the corresponding 
(Sobolev-)Hilbert spaces by 
\begin{equation*}
	H^s_1 \defn \set{f \in L_2(\mu) \colon 
    \norm[H_1^s]{f}^2 \,:=\, 
    \sum_{\ell \in \N_0} 
	\max\left\{1,\ell^{2s}\right\}\cdot |\inner{f}{b_\ell}_\nu|^2 <\infty}.
\end{equation*}
For $s\in\N$, the norm $\norm[H_1^s]{f}$ is equivalent to 
$\|f\|_{L_2(\mu)}+\|D^sf\|_{L_2(\mu)}$ for all the examples above, which justifies the name Sobolev space. 
However, note that these spaces consist solely of periodic functions for trigonometric polynomials $b_\ell=b_\ell^T$ and that we consider a weighted $L_2$-norm for Chebyshev polynomials $b_\ell=b_\ell^C$.

We now define the 
\emph{Sobolev spaces of dominating mixed smoothness} as the tensor product spaces 
$H_d^s:=H_1^s\otimes\cdots\otimes H_1^s$, 
which are all functions $f\in L_2([0,1]^d,\mu)$ such that the norm 
\[
\norm[H_d^s]{f}^2 \defn \sum_{\ell \in 
			\N_0^d} \, \left(\prod_{j=1}^d \max\left\{1,\ell_j^{2s}\right\}\right)\cdot |\inner{f}{b_\ell}_\nu|^2 
\]
is finite, where
\(b_\ell \defn \otimes_{j=1}^d b_{\ell_j}\) with \(\ell = (\ell_1, \dots, \ell_d)\in\N_0^d\).
By this definition, the functions $b_\ell$ 
are an orthonormal basis of $L_2([0,1]^d,\mu)$, 
with $\mu$ being the $d$-fold product measure 
of $\mu_1$, 
and an orthogonal basis of $H_d^s$.
This makes these spaces comparably simple to analyze, at least for $L_2$-approximation. 

In fact, 
denoting by $B_d^s$ the unit ball in $H_d^s$, 
it is
known that 
\begin{equation}\label{eq:best-Hs}
\sup_{f\in B_d^s}\,\inf_{g \in V_n} \norm[L_2(\mu)]{f-g} 
\;\lesssim\; n^{-s} (\log n)^{s(d-1)}, 
\end{equation}
where $V_n$ is the span of those $n$ functions $b_\ell$ that have the smallest norm in $H_d^s$, see e.g.~\cite{DTU18}.
The corresponding ``frequencies'' $\ell$ form a \emph{hyperbolic cross}, 
and one can choose 
$V_n=E(q,d)$ from Proposition~\ref{prop:space} in case of algebraic polynomials if $n=N(q,d)$, see~\eqref{eq:Nqd}.

From this we immediately see that, 
for some choice of sampling points and weights, a corresponding weighted LS method satisfies the same (optimal) bound whenever $s>1/2$. 
Even with suitable random points 
we obtain the bound $n^{-s} (\log n)^{sd}$, 
see~\cite[Coro.~1]{KU21}. 
(Here, $n$ is the number of points, unlike the random part of Theorem~\ref{thm:LS-L2}.)
For Smolyak's algorithm from Section~\ref{sec:SA} 
it is only known that the bound
$n^{-s} (\log n)^{(s+1)(d-1)}$ 
can be achieved in the case of Chebyshev polynomials, see~\cite[Remark~11]{BNR00}. 
In the case of trigonometric polynomials, 
it is known that a Smolyak-type algorithm satisfies 
the upper bound $n^{-s} (\log n)^{(s+1/2)(d-1)}$, but also that this bound cannot be improved, see~\cite[Section~5.3]{DTU18}. 
(We were not able to find these bounds also for algebraic polynomials, 
but we believe them to hold, too.)

Let us add that it was the prevalent conjecture that Smolyak's algorithm is an optimal method for (certain) tensor product Hilbert spaces such as $H_d^s$, see e.g.~Conjecture 
5.26 in \cite{DTU18}. 
Corollary~2 of~\cite{KU21} 
showed that this is not the case.

\medskip

Remarkably, the sampling densities $\rho_n$ needed in Theorem~\ref{thm:LS-L2} 
for the $V_n$ from~\eqref{eq:best-Hs}
can also be given 
(to some extent) 
explicitly. 
First, note that the trigonometric and the Chebyshev polynomials are uniformly bounded, 
i.e., $\sup_{\ell\in\N_0^d}\max\{\|b^T_\ell\|_\infty,\|b^C_\ell\|_\infty\}<\infty$.
This implies that one can simply choose 
$\rho_n=1$, see \cite[Remark~1]{KU21} or~\cite[Proposition~4.4]{SU23}. 
Hence, in this case, 
Theorem~\ref{thm:LS-L2} is about 
an unweighted LS method based on uniformly distributed points. 
However, for Chebyshev polynomials, the error is measured in a weighted $L_2$-norm, while it is the usual mean-squared error otherwise.
\\
The case of Legendre polynomials, i.e., 
orthogonal algebraic polynomials w.r.t.~the Lebesgue measure, 
has been treated in~\cite{CCMNT15,CDL13,CM17,GNZ20,NJZ17} 
in the context of \emph{randomized} least squares approximation based on random points.
It turned out that sampling w.r.t.~the 
inverse of the \emph{Christoffel function}, 
which is given by $
\max\{|g(x)|\colon g\in V_n, \|g\|_{L_2}=1\}$, 
seems to be a good choice, 
and that this function is equivalent to
the 
\emph{tensorized Chebyshev density}
\[
\widetilde{\rho}_n(t) 
\;:=\; \prod_{j=1}^d \frac{1}{\pi\sqrt{t_{j}(1-t_{j})}}
\]
for $t=(t_1,\dots,t_d)\in(0,1)^d$, see e.g.~\cite{GNZ20}. 
The density of Theorem~\ref{thm:LS-L2}, 
see also~\cite{KU21,Ull20}, 
appears to be 
different from $\widetilde{\rho}_n$ in general. 
(In fact, $\rho_n$ usually depends on $V_n$ and $F$.
) 
However, 
we will still use it in our experiments. 
It would be interesting to study when $\widetilde{\rho}_n$ can indeed be used in Theorem~\ref{thm:LS-L2}.

\medskip

Motivated by the above discussion, 
we will perform numerical experiments for the 
(unweighted/plain) least squares estimator 
\begin{equation} \label{eq:LS-uniform}
\LS{q,d} (f) \,:=\, 
\underset{g\in E(q,d)}{\rm argmin}\, \sum_{i=1}^{2\cdot N(q,d)} 
\vert f(x_i)-g(x_i) \vert^2
\end{equation}
with independent and uniformly distributed points $x_i\sim \cU[0,1]^d$, 
as well as the weighted least squares estimator
\begin{equation} \label{eq:LS-cheb}
\wLS{q,d} (f) \,:=\, 
\underset{g\in E(q,d)}{\rm argmin}\, \sum_{i=1}^{2\cdot N(q,d)} 
\left(\prod_{j=1}^d \sqrt{x_{i,j}(1-x_{i,j})}\right)\cdot\vert f(x_i)-g(x_i) \vert^2,
\end{equation}
where $x_i=(x_{i,1},\dots,x_{i,d})$ 
are chosen w.r.t.~the Chebyshev density~$\widetilde{\rho}_n$. 

\medskip

We will compare these two methods with Smolyak's algorithm from Section~\ref{sec:SA}. 
Recall that $\Sm{q,d}$ uses $N(q,d)=\dim(E(q,d))$ function evaluations, see Proposition~\ref{prop:space} and~\eqref{eq:Nqd}. 
Hence, $\LS{q,d}$ and $\wLS{q,d}$ use twice as many samples but \emph{project} onto the same space of polynomials~$E(q,d)$. 
We think this is a fair comparison, and did not find a reasonable alternative for using the same number of points. (We would need to take smaller spaces $V$ then.) 

Note that in this context, 
i.e., approximation by polynomials and iid points, 
we are only aware of quite bad error bounds for $\LS{q,d}$, 
and one may expect $\wLS{q,d}$ to perform better at least for $L_2$-approximation, 
see e.g.~\cite{CDL13,CM17}. 
We also do not know whether 
suitable LS methods 
achieve the (optimal) polynomial order $n^{-s}$ for uniform approximation in $F_d^s$ just like $\Sm{q,d}$, 
see Theorem~\ref{thm:SA} 
and~\cite{KPUU25-general} for a collection of known results.

Our experimental findings indicate, however, that 
these simple methods (even without the \emph{logarithmic oversampling}) are often superior to $\Sm{q,d}$, 
also for approximation in the uniform norm. 
Moreover, $\LS{q,d}$ appears to be better than $\wLS{q,d}$ in high dimensions. 
This is not yet supported by theory.

\bigskip

\begin{remark}[Discretization] \label{rem:disc} 
Theorem~\ref{thm:LS-Lp} is based on a 
certain type of \emph{discretization} of $L_p$-norms on $n$-dimensional spaces $V_n\subset C(D)$, 
see~\cite{KPUU25-uniform}. 
In fact, 
it follows from the 
existence of points $x_1,\dots,x_{4n}\in D$ 
such that, for all $f\in V_n$ and $p\ge2$, we have
$$\|f\|_{L_p(\mu)} \leq
83 \cdot  n^{\frac{1}{2}-\frac{1}{p}} \cdot \left(\frac{1}{4n} \sum_{k=1}^{4n} \abs{f(x_k)}^2 \right)^{1/2}.$$
This also implies 
that 
$\|f\|_{\infty} \;\lesssim\; \sqrt{n}\cdot \max_{k} \abs{f(x_k)}$, 
see~\cite[eq.~(2)]{KPUU25-uniform}.

Point sets $P_n\subset D$ 
such that $\|f\|_\infty \le C_n \max_{x\in P_n} |f(x)|$ 
for all $f\in V_n$ 
are called 
\textit{(weakly) admissible meshes}
or, if $C_n\asymp 1$, 
\textit{norming sets/meshes}. 
Especially for $V_n$ being multivariate polynomials, this attracted quite some attention, and we refer to~\cite{Bos18a,KKT23,Kroo11,NX23} for more on this and its relation to (least squares) approximation. 
Let us only highlight the recent work~\cite{DP23} where the existence of 
\emph{near-optimal} polynomial norming sets over arbitrary convex domains has been established. 
Explicit constructions of optimal norming sets seem to be unknown even for simple sets as the cube, see~\cite{Bos18b,BKSVZ24}. 
Moreover, note that such improved discretizations of the uniform norm do not lead to better upper bounds for least squares methods, but only for a (nonlinear and usually infeasible) \emph{least maximum estimator}, see~\cite[Remark~19]{KPUU25-uniform}. 
\end{remark}

\medskip

\section{Numerical results}\label{sec:experiments}

To test the theoretical findings discussed, 
we implement Smolyak's algorithm 
from Section~\ref{sec:SA}, as well as two variants of the 
(weighted) least squares algorithm 
from Section~\ref{sec:LS}.
We then compare them for various 
(randomly generated) test functions. 

For other numerical experiments for Smolyak's algorithm, we again refer to~\cite{BNR00,JMMV14}. 
For least squares based on random points, see e.g. \cite{BSU23,CCMNT15,CDL13,CM17,GNZ20,MNST14,NJZ17}. 
In those works, however, the performance is evaluated for an individual method and/or based on a theoretical worst-case error or the condition number of the corresponding linear system. 
We are not aware of a numerical comparison for specific functions.

See also the recent work~\cite{GHM25} which somehow combines the ideas from both settings.

\goodbreak

\subsection{Methodology}\label{subsec:method} 
We now describe the details of our numerical experiments. 
Let us first recall the details of the algorithms, 
which we will use under the names 

\goodbreak

\begin{itemize}[leftmargin=30mm] 
\itemsep=2mm
    \item[Smolyak: ] $\Sm{q,d}$, i.e., interpolation on a sparse grids based on Chebyshev nodes~\eqref{eq:chebyExtrema} which consist of $N(q,d)$ points, 
    \item[LS-uniform: ] $\LS{q,d}$, i.e., unweighted least squares from~\eqref{eq:LS-uniform} based on $2\cdot N(q,d)$ 
    uniformly distributed points, and 
    \item[LS-Chebyshev: ] $\wLS{q,d}$, i.e., weighted least squares from~\eqref{eq:LS-cheb} based on $2\cdot N(q,d)$ 
    random points w.r.t.~the Chebyshev density. 
\end{itemize}
For every $q$ and $d$, these random points are generated only once and then fixed for the rest of the simulation.
All 3 methods are exact on $E(q,d)$ from~Proposition~\ref{prop:space}, i.e., $A(g)=g$ for $g\in E(q,d)$, at least with 
probability one.
See Remark~\ref{rem:details} for details on the used soft- and hardware.

\medskip

To compare their performance for a given function 
$f\colon[0,1]^d\to\R$, we generate $M$ uniformly distributed points \(\{x_i\}_{i=1}^M \sim \cU[0,1]^{d}\) and compute 
\[
e_{\rm max}(A,f) \defn \max_{i\in\set{1,\ldots, M}} 
    \abs{f(x_i)-A(f)(x_i)} 
\]
and 
\[
e_{\rm mean}(A,f) \defn \left(\frac{1}{M}\sum_{i=1}^{M}
    \left(f(x_i)-A(f)(x_i)\right)^2\right)^{1/2}
\]
for $A\in\{\Sm{q,d},\LS{q,d},\wLS{q,d}\}$.
We choose $M=N(q,d)$ for $d\le 10$ and $M= 100\cdot N(d+2,d)$ for large $d$. 
We do not claim that this leads to a good estimation of the \emph{true error}, but all approximants are tested with the same points. 
\medskip

Since we want to study a worst-case setting, 
we apply the methods to several functions. 
In fact, we consider 
\begin{itemize}
\itemsep=2mm
    \item test functions from 9 different classes, see Section~\ref{sec:test},  
    \item with $Q=50$ random functions $g_1,\dots,g_Q$ from each class by choosing corresponding 
    parameters \(c, w \sim \cU[0,1]^d\) and rescale such that \(\norm[1]{c}=d\). 
    \item All random points in the algorithms, as well as in $e_{\rm max}$ and $e_{\rm mean}$, stay fixed.
\end{itemize}
We then compute 
\[
    e_{\rm *}^{\rm wc}(A) \defn \max_{j\in\set{1,\ldots, Q}}  e_{\rm *}(A,g_j)
\]
for $A\in\{\Sm{q,d},\LS{q,d},\wLS{q,d}\}$ 
and ${\rm *}\in\{{\rm max},{\rm mean}\}$.

\medskip
\goodbreak

In the following section, 
we will report on our numerical results 
for the individual classes of functions and for several choices of $q$ and $d$. In particular, we  
\begin{itemize}
\itemsep=2mm
    \item denote the number $k:=q-d$ as \emph{scale}, and consider 
    \item scale up to $k=9$ for $d\le10$ in Section~\ref{sec:num-low}, 
        as well as 
    \item dimension up to $d=100$ for scale $k=1,2$ in Section~\ref{sec:num-high}, and
    \item \emph{stability} when data is given by random noise in Section~\ref{sec:num-noise}.
\end{itemize}

To reproduce our findings, we refer the reader to our ``plug-and-play'' 
implementation available on GitHub.
See also Section~\ref{sec:discussion} for running times and possible drawbacks. 

\medskip

\begin{remark}[Details on Hardware and Software] \label{rem:details}
All experiments were carried out on a high performance computing (HPC) cluster featuring dual-socket AMD EPYC 7742 processors, each operating at 2.25\,GHz, and equipped with 1\,TB of RAM.
Throughout the project, we use \textit{Numpy} and \textit{Scipy} \cite{HMW20, VGO20}. 
For finding the least-squares solution we employ
the \textit{gelsy} driver from the \textit{Lapack} library \cite{ABB99}.
For Smolyak's algorithm, we employ the Tasmanian library \cite{Sto15,SLBM13}. 
The spaces $E(q,d)$ are implemented with the Chebyshev basis based on \cite{CS13}. 
\end{remark}

\bigskip
\goodbreak
\newpage

\subsection{Test functions} \label{sec:test}

We use the following families of test functions that are based on \cite{Gen84,Gen87,SB13}, each defined on the $d$-dimensional unit-cube $[0,1]^d$, and specified by parameters $c,w\in[0,1]^d$.


\allowdisplaybreaks[4]
\begin{align*}
    \begin{array}{l l}
        \text{1. Bimodal Gaussian:} & f_{1}(x) = 
        \phi\left(x-w\right)+\phi\left(x+w-1\right) 
        \\[12pt]
        & \quad\text{with }\quad \phi(x) := 
        \exp \left(-\frac{50}{d}\sum_{i=1}^{d}\left(c_ix_i\right)^2\right)
        \\[20pt]
        \text{2. Continuous:} & f_2(x) = \exp\left( - \frac{1}{d} 
        \sum_{i=1}^{d} c_i \abs{x_i - w_i} \right) \\[20pt]
        \text{3. Corner Peak:} & f_3(x) = \left( 1 + 
        \sum_{i=1}^{d} c_i x_i 
        \right)^{-(d+1)} \\[20pt]
        \text{4. Discontinuous:} & f_4(x) = 
        \begin{cases}
            0, & x_1 > w_1 \lor x_2 > w_2, \\
            \exp\left( \sum_{i=1}^{d} c_i x_i \right), & \text{else}
        \end{cases} \\[24pt]
        \text{5. Gaussian:} & f_5(x) = \phi(x-w) \\[20pt]
        \text{6. Geometric Mean:} & f_6(x) = 
        \left(1+1/d\right)^d  
        \prod_{i=1}^{d}\left(c_ix_i+w_i\right)^{1/d}
        \\[20pt]
        \text{7. Oscillatory:} & f_7(x) = \cos\left( 2\pi w_1 + 
        3\sum_{i=1}^{d} c_i x_i \right) \\[20pt]
        \text{8. Product Peak:} & f_8(x) = \prod_{i=1}^{d} 
        \left( c_i^{-2} + (x_i - w_i)^2 \right)^{-1} \\[20pt]
        \text{9. Ridge Product:} & f_9(x) = \prod_{i=1}^{d}
        \frac{\abs{4x_i-2-w_i}+c_i}{1+c_i}
    \end{array}
\end{align*}

\begin{figure}[H]
	\centering
	\includegraphics[width=\linewidth]{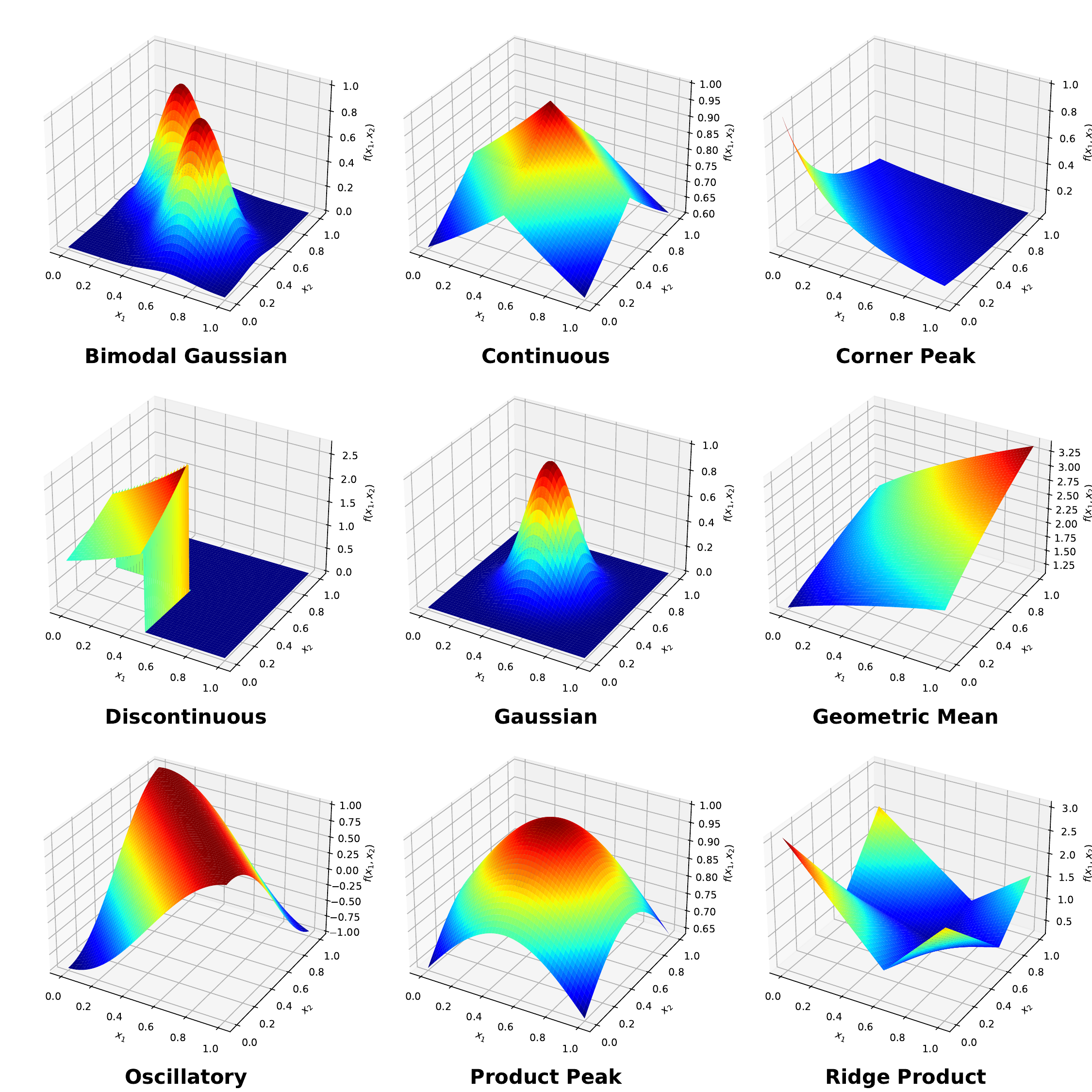}
	\caption{Visualization of the functions for $d=2$.}
    \label{fig:function_visualization}
\end{figure}

Note that $f_1$ (Bimodal Gaussian), $f_6$ (Geometric Mean) and $f_9$ (Ridge 
Product) are called \emph{Zhou}, \emph{Morokoff Calfisch 1} and 
\emph{G-Function} in \cite{SB13} and in our code. 
We adapted them slightly to allow for a parametric choice.
The function classes $f_2, f_3, f_4, f_5, f_7$ and $f_8$ are also known as the \emph{Genz 
Integrand Families} and were introduced by Genz in \cite{Gen84, Gen87}. 
Let us add that we also tried other functions of~\cite{SB13} but without any interesting outcome.
See Figure~\ref{fig:function_visualization} for prototypical examples of our test functions.

\subsection{Approximation in low dimensions} 
\label{sec:num-low}

We now present some of our results for 
dimensions
up to 
$d=10$ with ``large'' scales up to 
9. 
Hence, the maximal numbers of points used 
are given, e.g., by $N(12,3)=13953$, $N(13,5)=51713$ and $N(16,10)=171425$ for Smolyak, while both LS use twice as many.
First, we present a table giving the precise values of $e_{\rm max}^{\rm wc}$ and $e_{\rm mean}^{\rm wc}$ for all algorithms, functions and scales in dimension 10. 
(Tables for all smaller dimensions can be found on GitHub.) 
Afterwards, we give some illustrations by plotting some 
individual results, i.e., rows from the tables.

\begin{table}[H]
	\centering
    \makebox[\linewidth]{%
    	\begin{adjustbox}{width=\linewidth}
\begin{tabular}{ll|rr|rr|rr|rr|}
 &   \multicolumn{1}{c}{} & \multicolumn{2}{c}{Scale 3} & \multicolumn{2}{c}{Scale 4} & \multicolumn{2}{c}{Scale 5} & \multicolumn{2}{c}{Scale 6}\\
 &   \multicolumn{1}{c}{} & \multicolumn{1}{c}{$e_{\rm max}^{\rm wc}$} & \multicolumn{1}{c}{$e_{\rm mean}^{\rm wc}$} & \multicolumn{1}{c}{$e_{\rm max}^{\rm wc}$} & \multicolumn{1}{c}{$e_{\rm mean}^{\rm wc}$} & \multicolumn{1}{c}{$e_{\rm max}^{\rm wc}$} & \multicolumn{1}{c}{$e_{\rm mean}^{\rm wc}$} & \multicolumn{1}{c}{$e_{\rm max}^{\rm wc}$} & \multicolumn{1}{c}{$e_{\rm mean}^{\rm wc}$}\\
\toprule
\multirow{3}{*}{\thead[l]{\tiny\textbf{Bim. Gauss.}\\$Q=50$}} & Smolyak & 1.21e+00 & 2.73e-01  & 9.42e-01 & 1.45e-01  & 8.41e-01 & 6.38e-02  & 4.33e-01 & 4.68e-02\\
 & LS-Uniform & \first{5.24e-01} & \first{6.01e-02}  & \first{4.85e-01} & \first{4.17e-02}  & 1.13e+00 & 2.66e-02  & 1.52e+00 & 1.71e-02\\
 & LS-Chebyshev & 6.88e-01 & 6.53e-02  & 5.97e-01 & 4.47e-02  & \first{4.92e-01} & \first{2.57e-02}  & \first{3.66e-01} & \first{1.55e-02}\\
\midrule
\multirow{3}{*}{\thead[l]{\tiny\textbf{Cont.}\\$Q=50$}} & Smolyak & 2.33e-02 & 4.16e-03  & \first{9.30e-03} & 1.55e-03  & \first{5.37e-03} & 6.20e-04  & 2.79e-03 & \first{1.92e-04}\\
 & LS-Uniform & 2.36e-02 & \first{3.73e-03}  & 1.20e-02 & 1.48e-03  & 1.07e-02 & 6.01e-04  & 6.48e-03 & 2.37e-04\\
 & LS-Chebyshev & \first{2.23e-02} & 4.08e-03  & 1.03e-02 & \first{1.45e-03}  & 5.61e-03 & \first{5.39e-04}  & \first{2.34e-03} & 1.95e-04\\
\midrule
\multirow{3}{*}{\thead[l]{\tiny\textbf{Corner Peak}\\$Q=50$}} & Smolyak & 4.08e-04 & \first{1.11e-05}  & 5.66e-04 & \first{1.08e-05}  & 1.16e-02 & 1.67e-04  & 2.88e-02 & 4.11e-04\\
 & LS-Uniform & \first{3.87e-04} & 2.09e-05  & \first{5.40e-04} & 1.26e-05  & \first{9.89e-04} & \first{9.84e-06}  & \first{2.38e-03} & \first{8.90e-06}\\
 & LS-Chebyshev & 1.62e-03 & 2.92e-04  & 1.33e-03 & 2.53e-04  & 4.98e-03 & 3.72e-04  & 6.19e-03 & 3.31e-04\\
\midrule
\multirow{3}{*}{\thead[l]{\tiny\textbf{Discont.}\\$Q=50$}} & Smolyak & 4.14e+03 & 2.24e+02  & 4.14e+03 & \first{1.37e+02}  & 2.92e+03 & \first{1.00e+02}  & 5.01e+03 & \first{7.81e+01}\\
 & LS-Uniform & \first{3.36e+03} & \first{1.93e+02}  & \first{2.89e+03} & 1.43e+02  & 2.62e+03 & 1.13e+02  & 4.37e+03 & 9.03e+01\\
 & LS-Chebyshev & 3.79e+03 & 4.90e+02  & 3.60e+03 & 3.37e+02  & \first{2.60e+03} & 2.30e+02  & \first{3.39e+03} & 1.75e+02\\
\midrule
\multirow{3}{*}{\thead[l]{\tiny\textbf{Gauss.}\\$Q=50$}} & Smolyak & 7.29e-01 & 8.59e-02  & 6.10e-01 & 5.78e-02  & 5.69e-01 & 3.47e-02  & 4.24e-01 & 2.04e-02\\
 & LS-Uniform & \first{5.32e-01} & \first{4.63e-02}  & \first{5.28e-01} & \first{2.98e-02}  & 1.01e+00 & \first{1.82e-02}  & 9.98e-01 & 1.08e-02\\
 & LS-Chebyshev & 7.05e-01 & 5.57e-02  & 5.66e-01 & 3.15e-02  & \first{4.67e-01} & 1.96e-02  & \first{3.04e-01} & \first{1.07e-02}\\
\midrule
\multirow{3}{*}{\thead[l]{\tiny\textbf{Geo. Mean}\\$Q=50$}} & Smolyak & 2.11e-01 & 3.30e-02  & 1.24e-01 & 1.11e-02  & 5.30e-02 & 3.61e-03  & 1.45e-02 & 1.10e-03\\
 & LS-Uniform & \first{1.33e-01} & \first{1.46e-02}  & 7.36e-02 & \first{4.07e-03}  & 2.87e-02 & \first{8.66e-04}  & 1.27e-02 & \first{1.82e-04}\\
 & LS-Chebyshev & 1.70e-01 & 3.77e-02  & \first{5.84e-02} & 1.26e-02  & \first{1.60e-02} & 2.92e-03  & \first{7.23e-03} & 8.24e-04\\
\midrule
\multirow{3}{*}{\thead[l]{\tiny\textbf{NOISE}\\$Q=50$}} & Smolyak & 1.06e-05 & 1.74e-06  & 2.09e-05 & 2.14e-06  & 3.98e-05 & 3.22e-06  & 6.21e-05 & 4.71e-06\\
 & LS-Uniform & 7.99e-07 & \first{1.25e-07}  & 1.42e-06 & \first{1.31e-07}  & 3.67e-06 & \first{1.47e-07}  & 5.65e-06 & \first{1.65e-07}\\
 & LS-Chebyshev & \first{7.85e-07} & 1.77e-07  & \first{1.01e-06} & 1.76e-07  & \first{1.15e-06} & 1.70e-07  & \first{1.29e-06} & 1.68e-07\\
\midrule
\multirow{3}{*}{\thead[l]{\tiny\textbf{Osci.}\\$Q=50$}} & Smolyak & 4.10e+01 & 3.10e+00  & 4.89e+01 & 2.10e+00  & 4.38e+01 & 1.44e+00  & 4.67e+01 & 9.71e-01\\
 & LS-Uniform & 6.72e+00 & \first{1.12e+00}  & 1.70e+01 & \first{1.13e+00}  & 4.13e+01 & \first{1.09e+00}  & 6.83e+01 & \first{8.65e-01}\\
 & LS-Chebyshev & \first{6.05e+00} & 1.46e+00  & \first{7.55e+00} & 1.44e+00  & \first{7.39e+00} & 1.34e+00  & \first{8.18e+00} & 1.20e+00\\
\midrule
\multirow{3}{*}{\thead[l]{\tiny\textbf{Prod. Peak}\\$Q=50$}} & Smolyak & \first{1.53e-02} & 1.89e-03  & \first{4.83e-03} & 4.05e-04  & \first{1.79e-03} & 6.80e-05  & \first{2.96e-04} & 1.34e-05\\
 & LS-Uniform & 3.22e-02 & \first{1.41e-03}  & 8.72e-03 & \first{2.71e-04}  & 3.74e-03 & \first{5.71e-05}  & 1.61e-03 & \first{1.08e-05}\\
 & LS-Chebyshev & 1.99e-02 & 2.59e-03  & 8.62e-03 & 5.31e-04  & 2.12e-03 & 9.73e-05  & 5.25e-04 & 1.85e-05\\
\midrule
\multirow{3}{*}{\thead[l]{\tiny\textbf{Ridge Prod.}\\$Q=50$}} & Smolyak & 4.96e+01 & 2.83e+00  & 6.16e+01 & 2.16e+00  & 4.62e+01 & 1.64e+00  & 5.71e+01 & 1.08e+00\\
 & LS-Uniform & \first{1.98e+01} & \first{1.64e+00}  & \first{1.61e+01} & \first{9.34e-01}  & 1.40e+01 & \first{5.31e-01}  & 1.15e+01 & \first{2.94e-01}\\
 & LS-Chebyshev & 2.68e+01 & 6.33e+00  & 1.98e+01 & 3.37e+00  & \first{1.00e+01} & 1.63e+00  & \first{7.43e+00} & 7.75e-01\\
\bottomrule
\end{tabular}

    	\end{adjustbox}%
    }
	\caption{The errors $e_{\rm max}^{\rm wc}$ and $e_{\rm mean}^{\rm wc}$ 
    for $d=10$ and 
    each scale 
    for the different algorithms. 
    Lowest value in bold.}
	\label{tab:dim10_results}
\end{table}

\bigskip
\goodbreak

Table~\ref{tab:dim10_results} shows that 
one of the least squares methods often produces the smallest error.

The following Figures~\ref{fig:figures_dim10_1} and \ref{fig:figures_dim10_2} visualize all the rows from Table~\ref{tab:dim10_results}, 
except for ``continuous'', which did not show anything interesting, 
and ``noise'', which will be discussed in Section~\ref{sec:num-noise}.
Note that all three algorithms show comparable convergence in all the cases considered. 
Additional interesting phenomena could be that 
\begin{itemize}
\itemsep=1mm
    \item LS-uniform often has smaller values of $e_{\rm mean}^{\rm wc}$ compared to LS-Chebyshev, and vice versa for $e_{\rm max}^{\rm wc}$. 
    \item for some ``non-smooth'' functions (like ``Ridge Product'' and ``Corner Peak'') it seems that LS has an asymptotic advantage.
    \item for (bimodal) Gaussian, LS-uniform seems to have problems for large scales.
\end{itemize}
The latter is even more clear in  
lower dimensional cases that will be discussed afterwards.
Note that the results for SA for $d=10$ are not identical to~\cite{BNR00} due to a different normalization. 
We do not believe this is essential, as our results should be understood as comparisons.

\begin{figure}[H]
	\centering
	\includegraphics[width=\widthoffigures\linewidth]{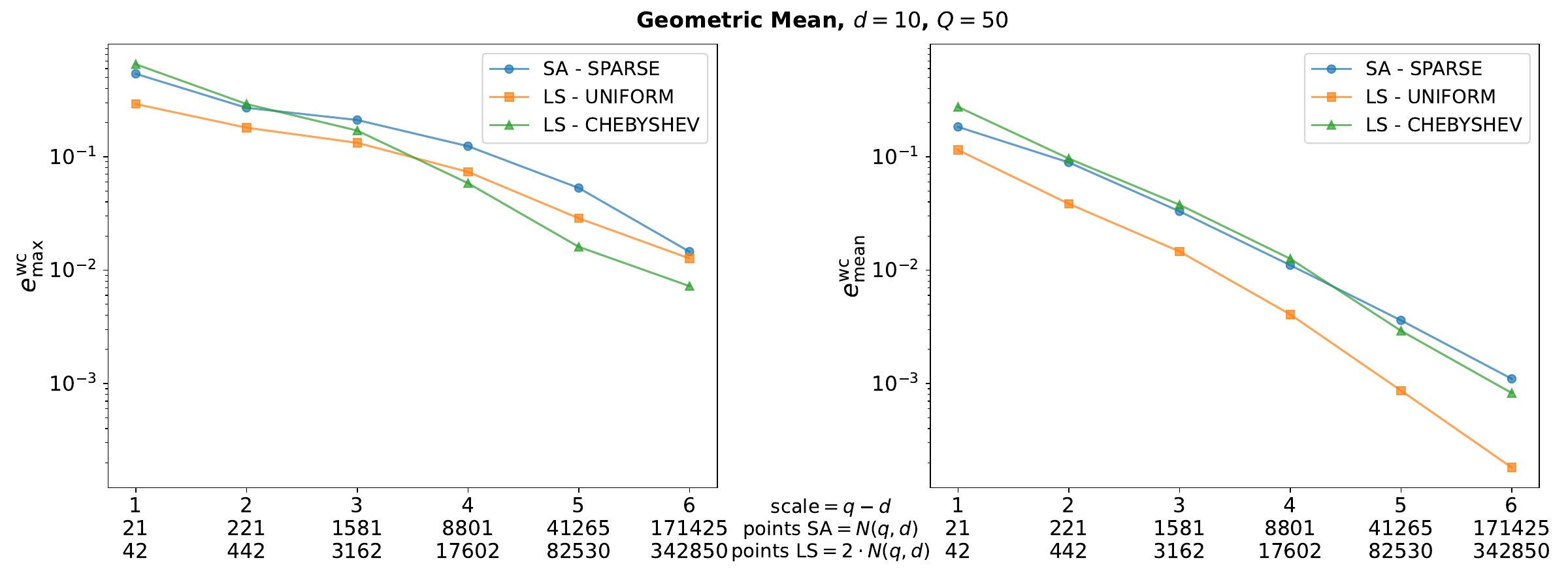}
	\includegraphics[width=\widthoffigures\linewidth]{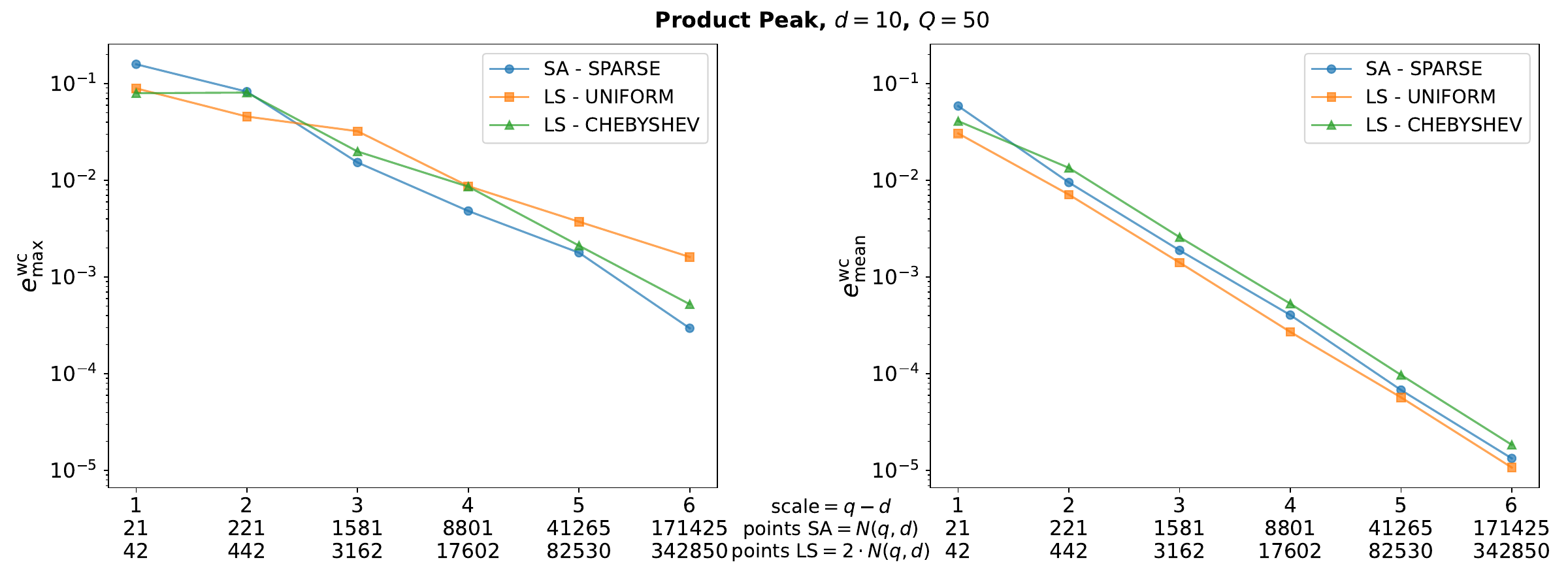}
	\includegraphics[width=\widthoffigures\linewidth]{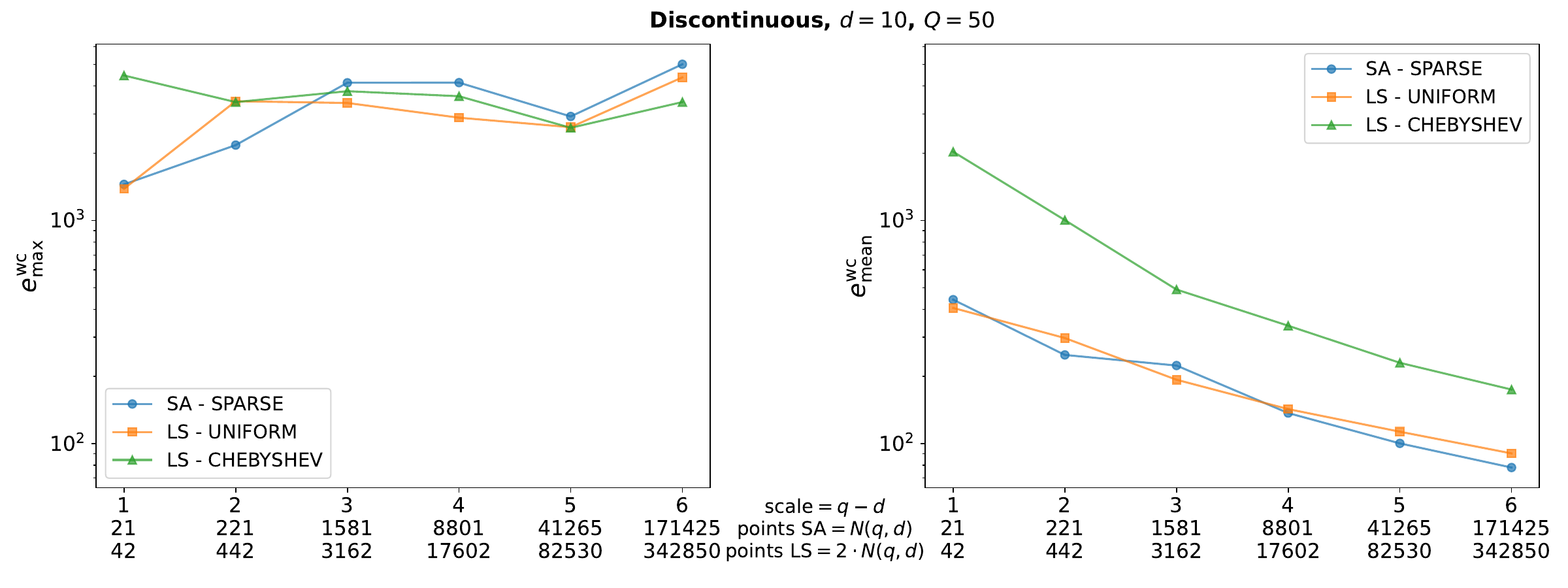}
	\includegraphics[width=\widthoffigures\linewidth]{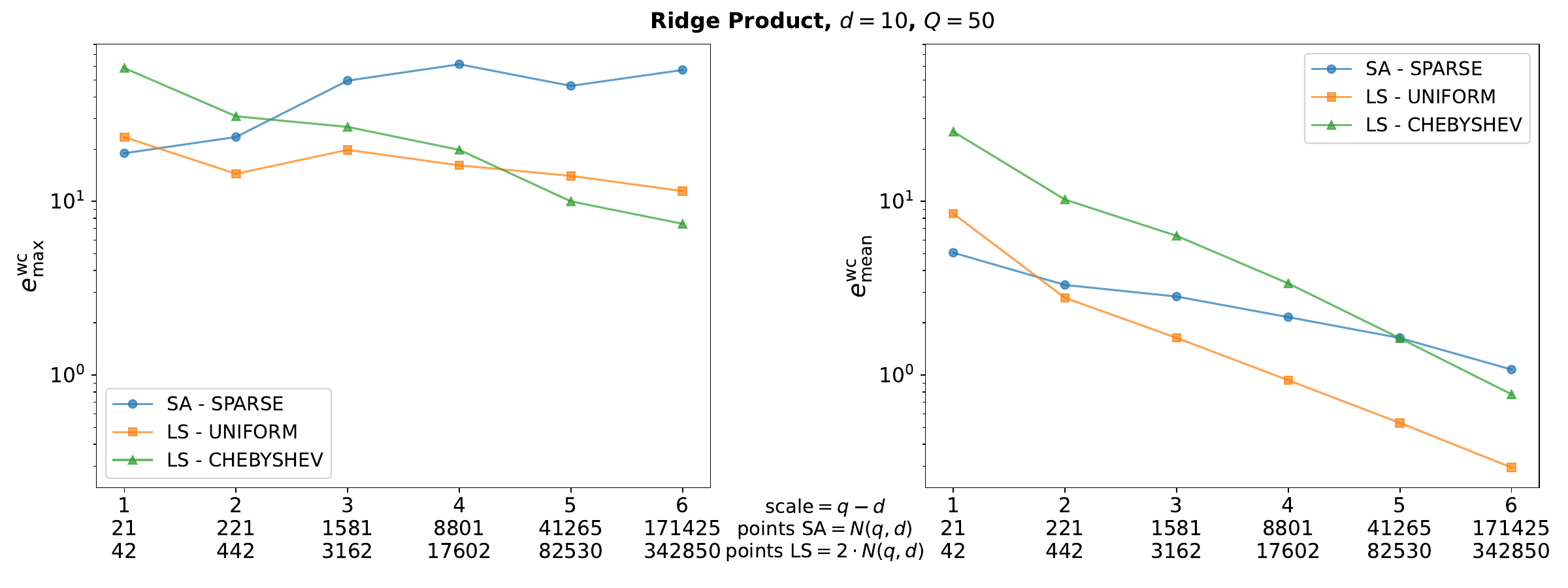}

	\caption{Some functions at fixed dimension \(d=10\).}
    \label{fig:figures_dim10_1}
\end{figure}

\begin{figure}[H]
	\centering
	\includegraphics[width=\widthoffigures\linewidth]{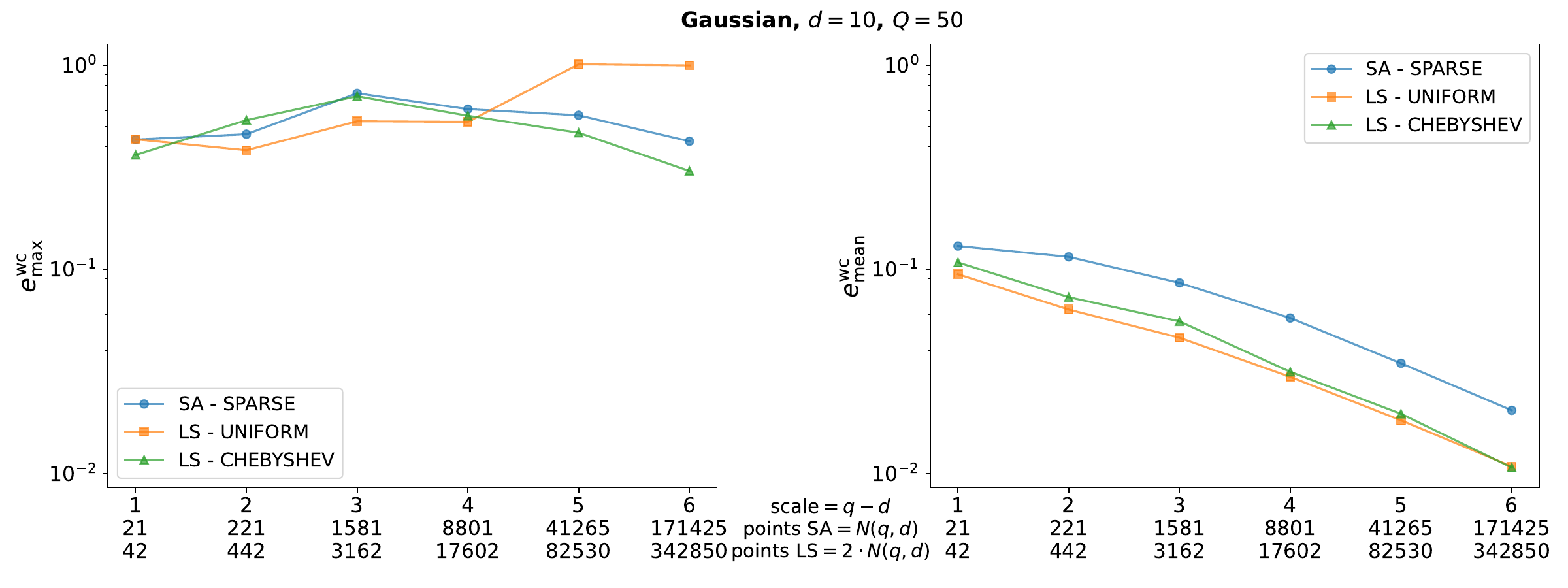}
	\includegraphics[width=\widthoffigures\linewidth]{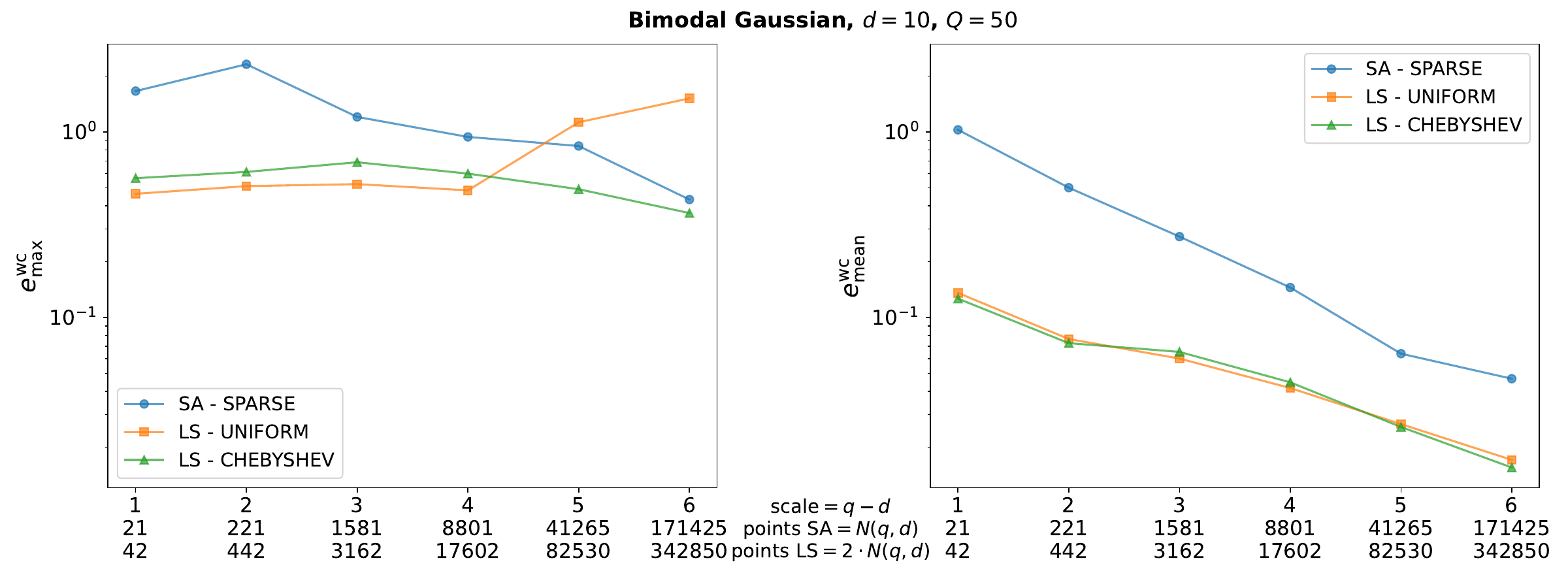}
	\includegraphics[width=\widthoffigures\linewidth]{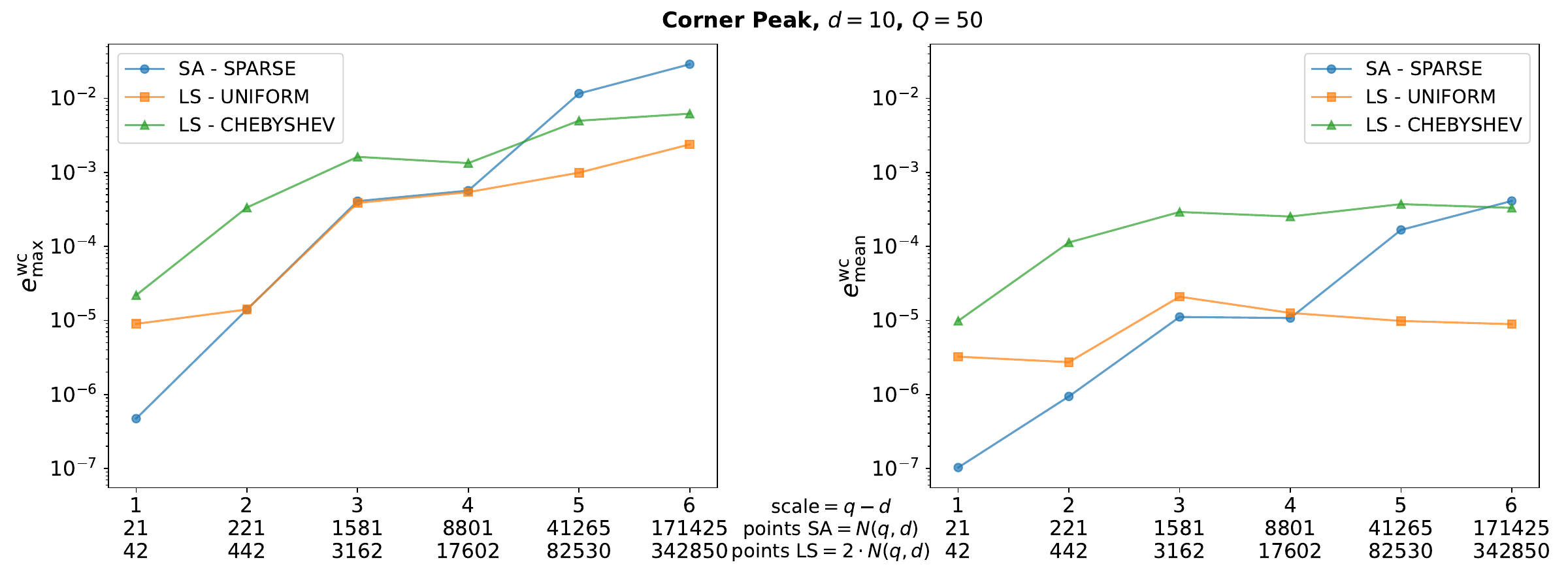}
    \includegraphics[width=\widthoffigures\linewidth]{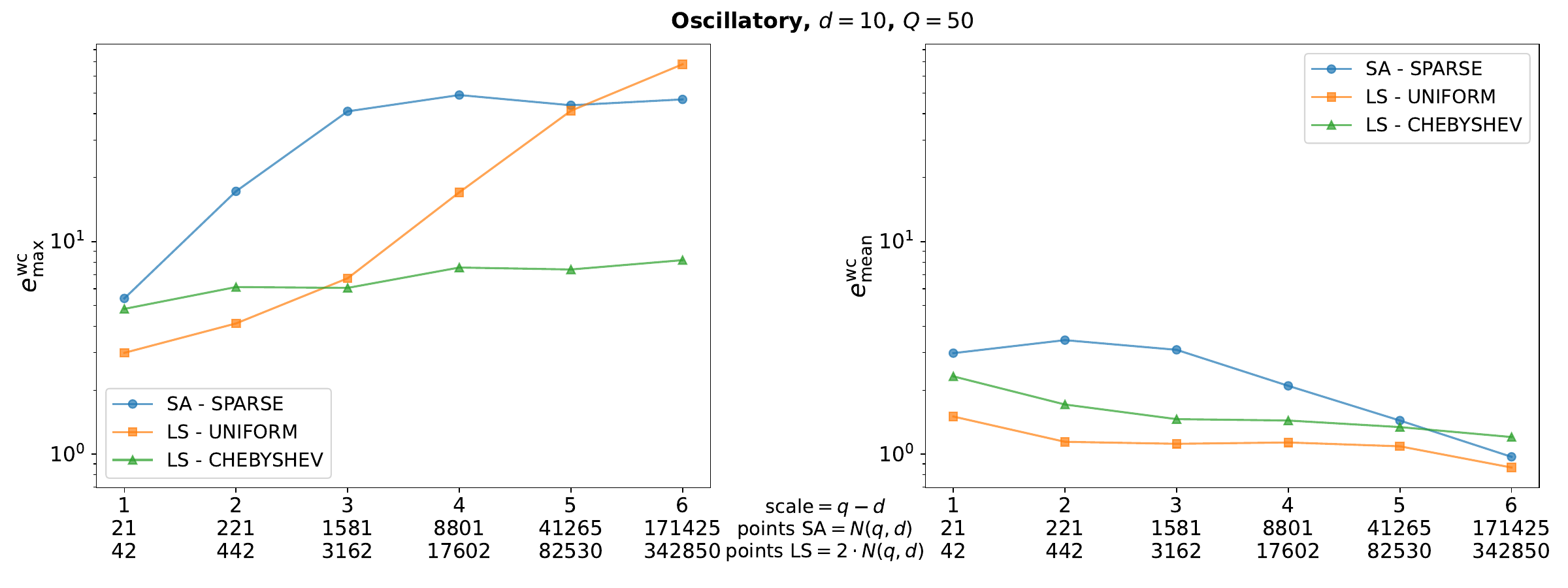}
    
    \caption{Some functions at fixed dimension \(d=10\).}
    \label{fig:figures_dim10_2}
\end{figure}

We now show some examples in smaller dimensions 
($d=3$ and $d=5$), 
and therefore larger scales, 
that appear to be interesting. 
(Figures for all other $d$ and $q$ can be found on GitHub.)

\begin{figure}[H]
	\centering
	\includegraphics[width=\widthoffigures\linewidth]{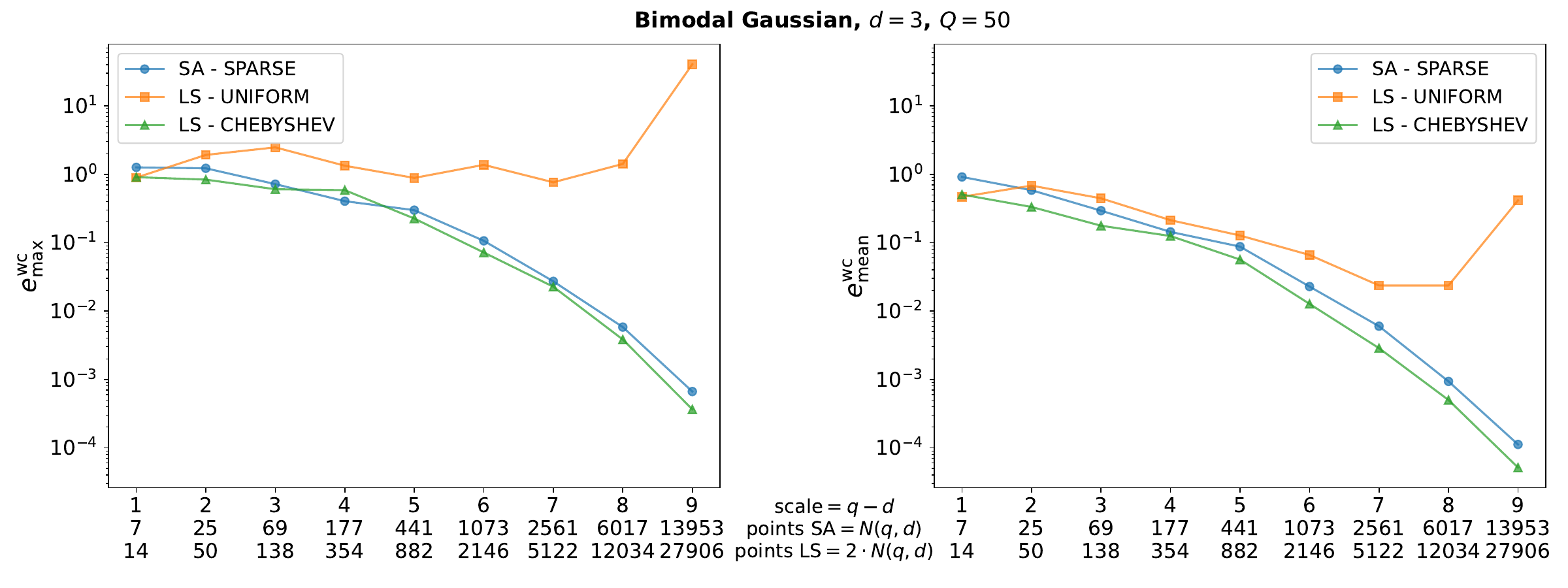}
    \includegraphics[width=\widthoffigures\linewidth]{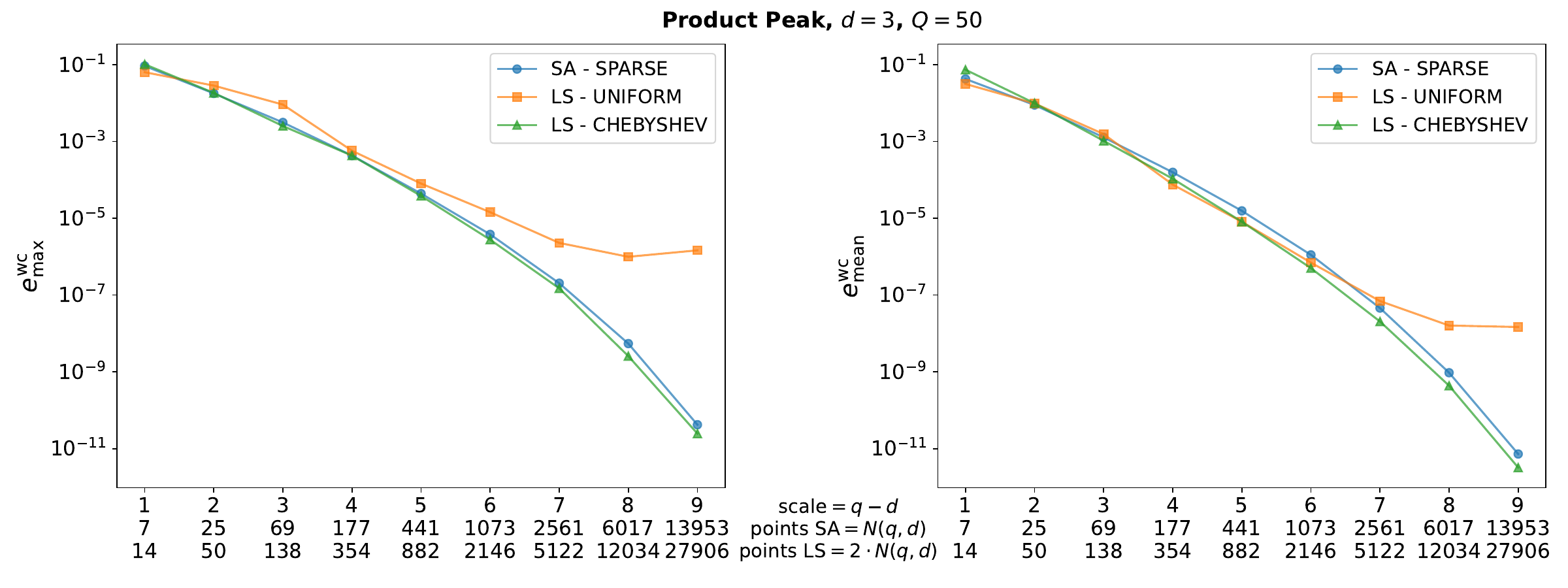}
    \includegraphics[width=\widthoffigures\linewidth]{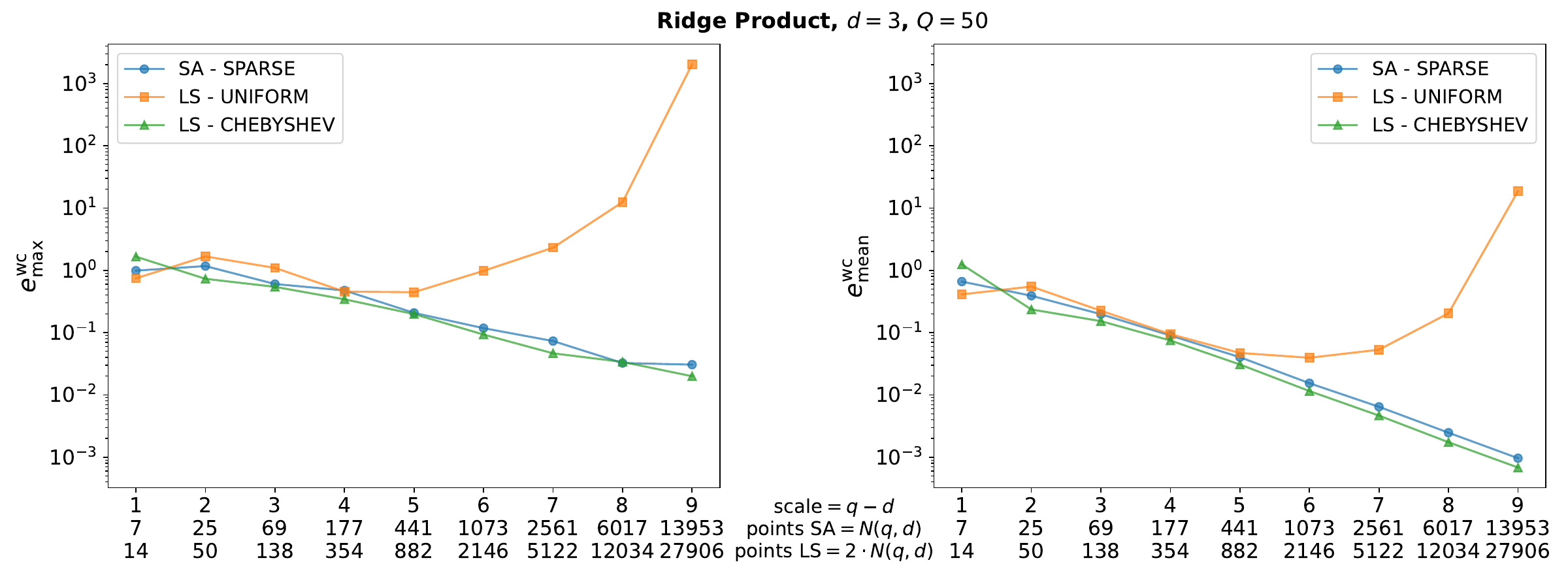}

\vspace{-4mm}

	\caption{Some functions at fixed dimension \(d=3\).}
	\label{fig:scale_1_9_dim3}
\end{figure}

We see that, again, SA and LS-Chebyshev show a comparable behavior, 
with LS-Chebyshev usually having a smaller error.
However, we also see that LS-uniform is worse
(except for exceptional cases such as ``Corner Peak'').
The difference is sometimes quite large, 
see Figure~\ref{fig:scale_1_9_dim3}, 
and it seems that, for large scales, 
LS-uniform would require a larger number of points to catch up. 
This is a known phenomenon, see e.g.~ \cite{CM17,MNST14}.

\begin{figure}[H]
	\centering
	\includegraphics[width=\widthoffigures\linewidth]{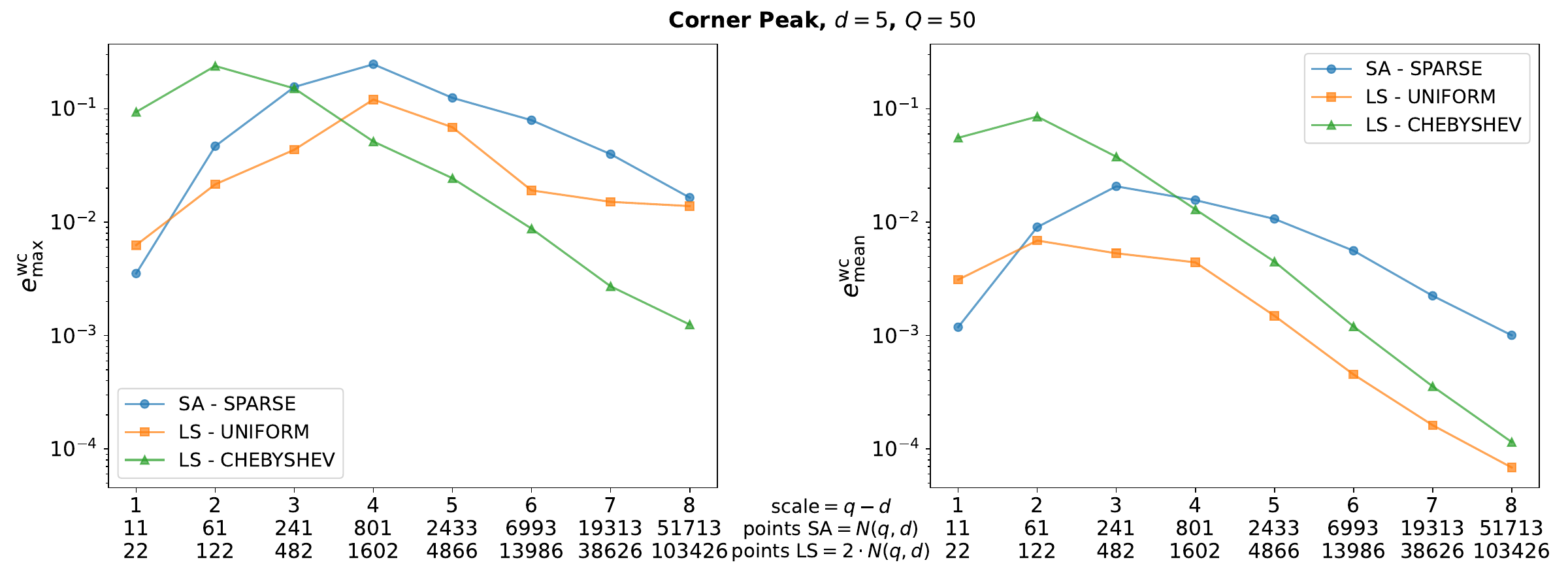}
    \includegraphics[width=\widthoffigures\linewidth]{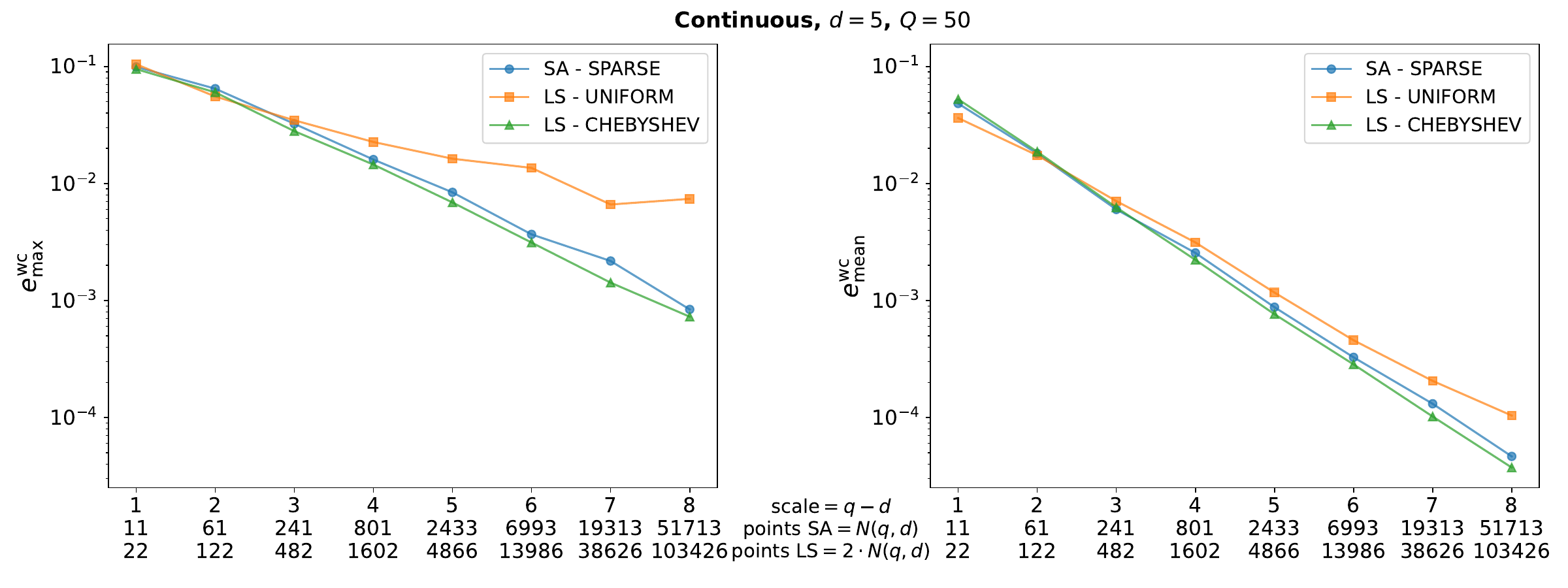}
    \includegraphics[width=\widthoffigures\linewidth]{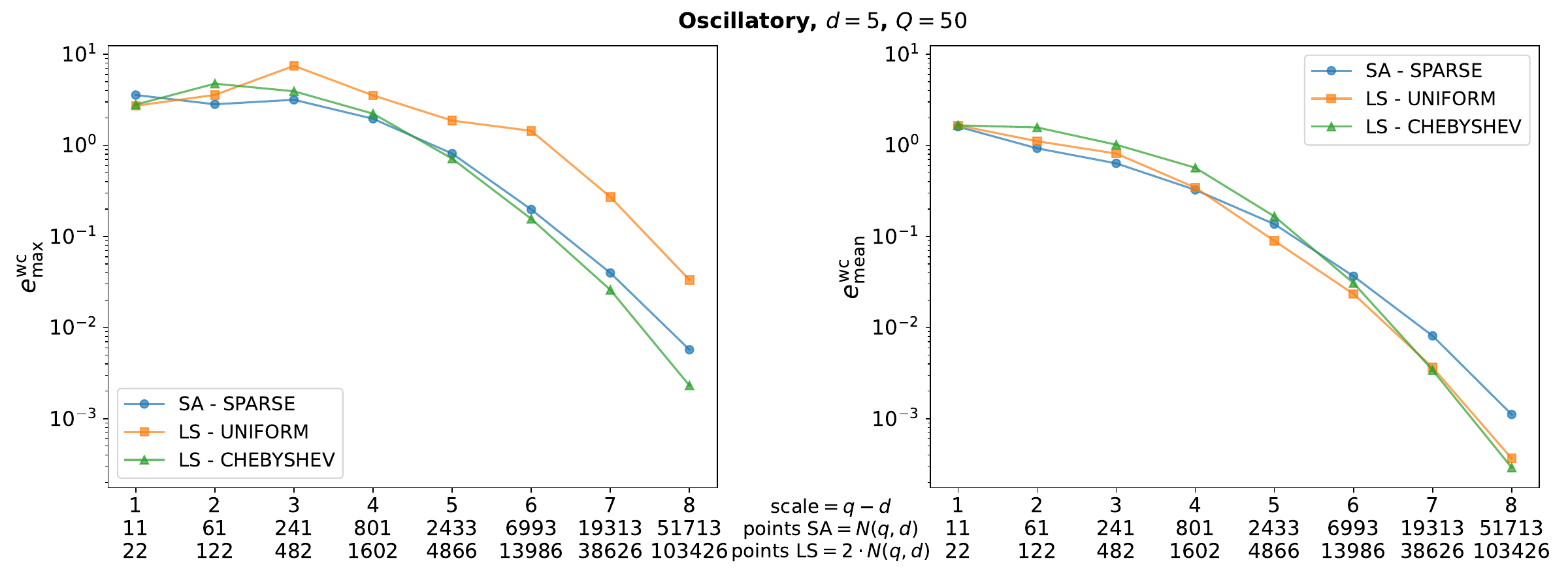}

	\caption{Some functions at fixed dimension \(d=5\).}
	\label{fig:scale_1_8_dim5}
\end{figure}

\goodbreak

\vspace{-3mm}

\subsection{Approximation in high dimensions} 
\label{sec:num-high}

\begin{figure}[H]
	\centering
	\includegraphics[width=\widthoffigures\linewidth]{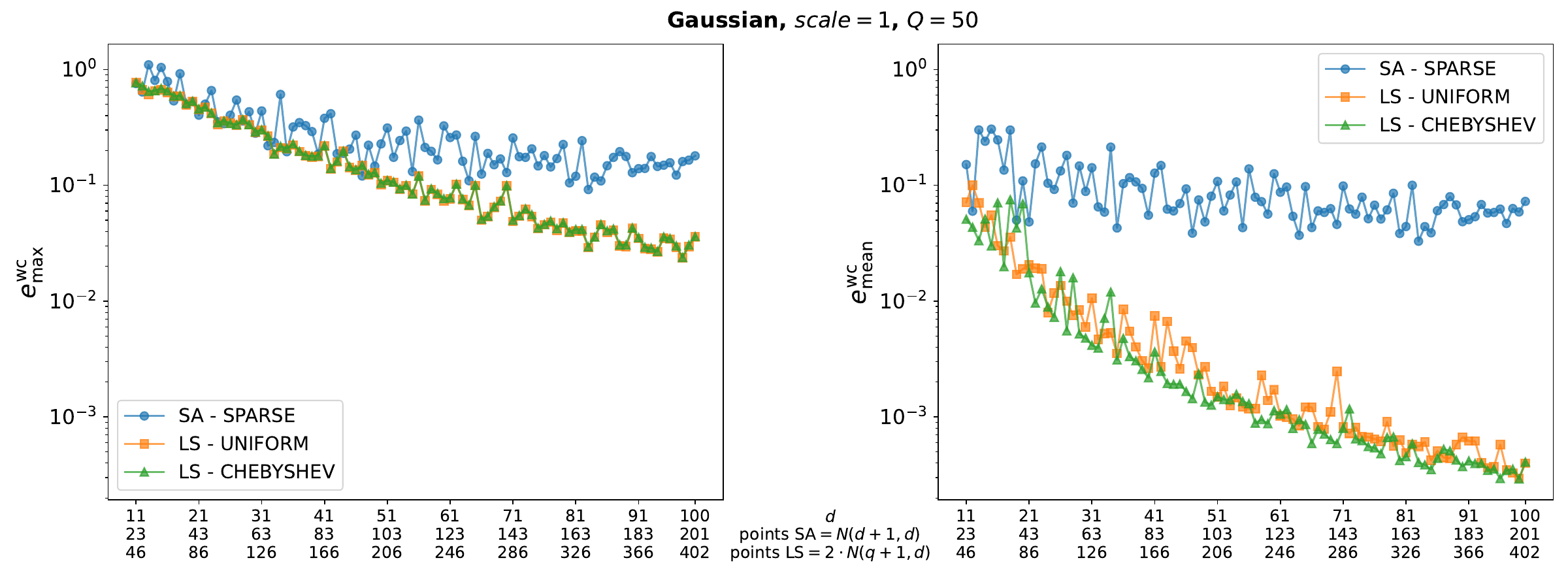}
	\includegraphics[width=\widthoffigures\linewidth]{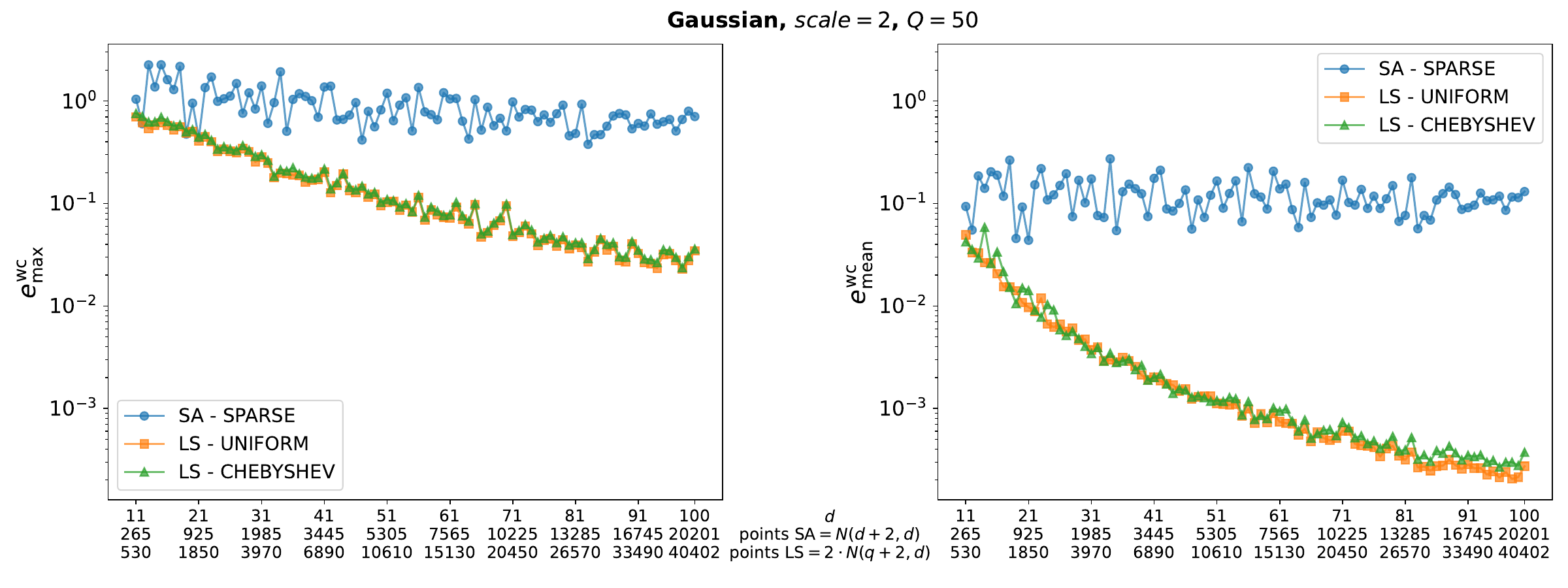}
	\includegraphics[width=\widthoffigures\linewidth]{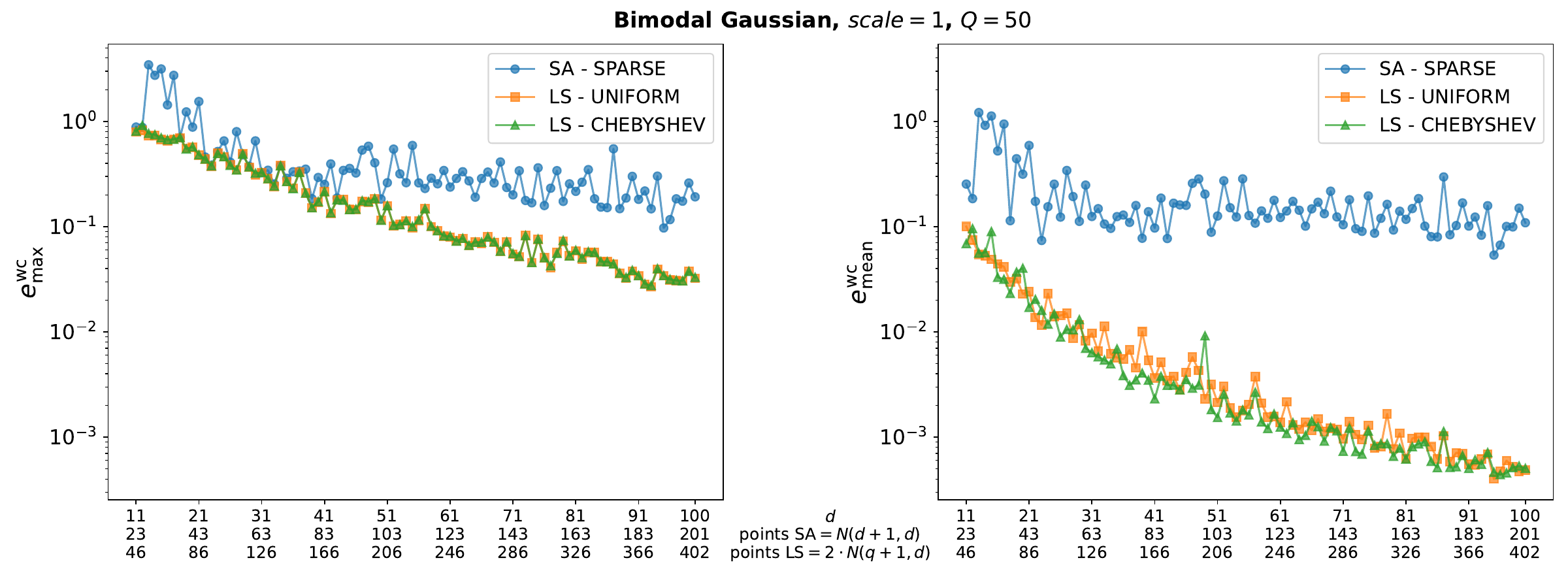}
	\includegraphics[width=\widthoffigures\linewidth]{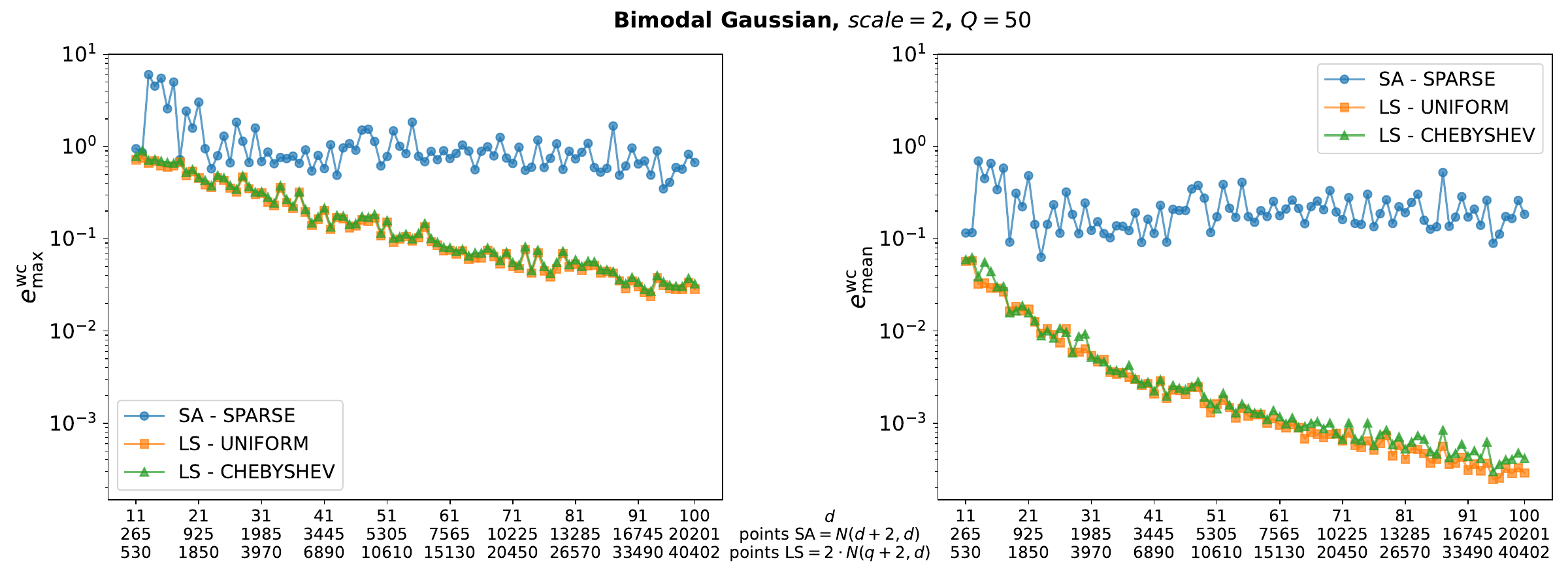}

    \caption{Some functions for fixed scale 1 and 2.}
	\label{fig:scale12_highdim_some_f_1}
\end{figure}

\begin{figure}[H]
	\centering

    \includegraphics[width=\widthoffigures\linewidth]{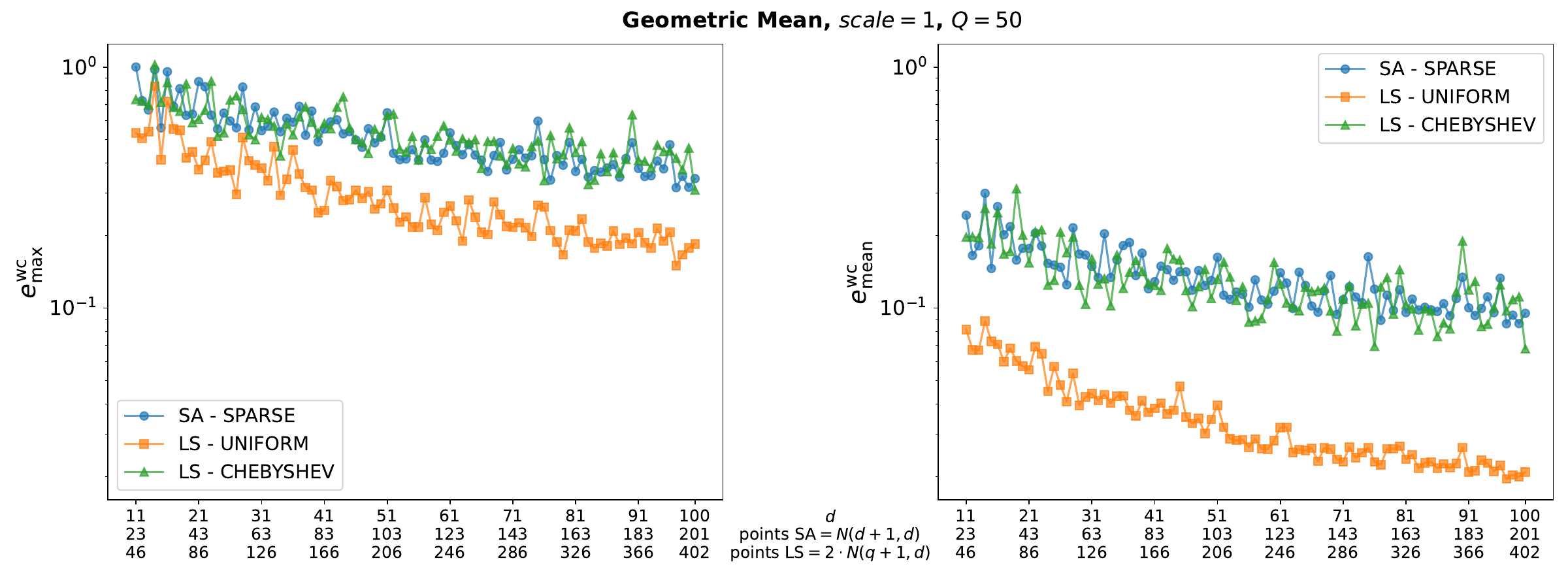}
	\includegraphics[width=\widthoffigures\linewidth]{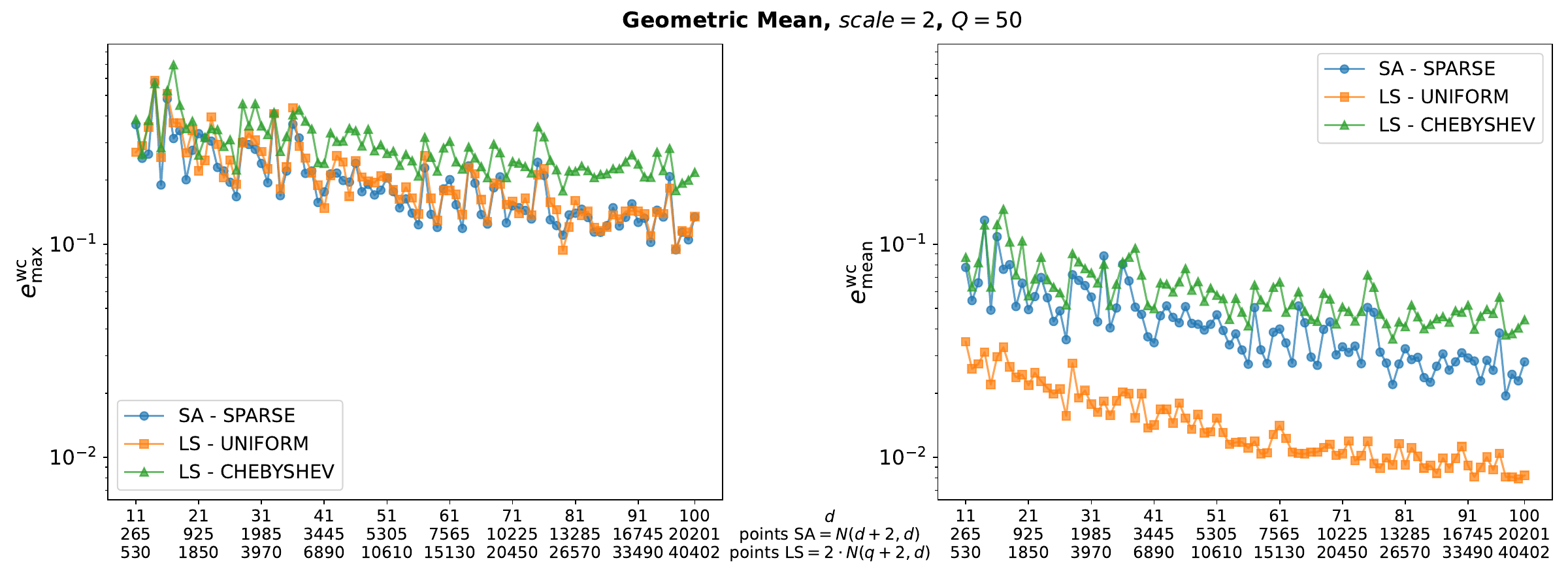}
    \includegraphics[width=\widthoffigures\linewidth]{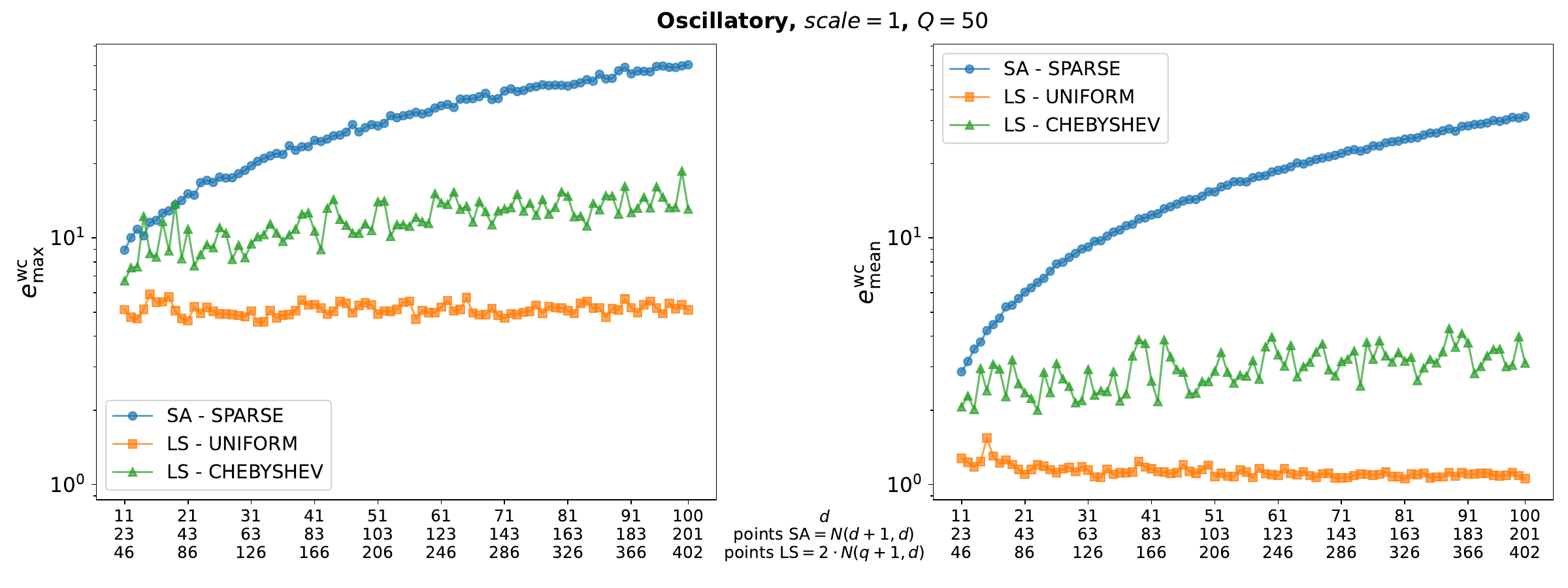}
	\includegraphics[width=\widthoffigures\linewidth]{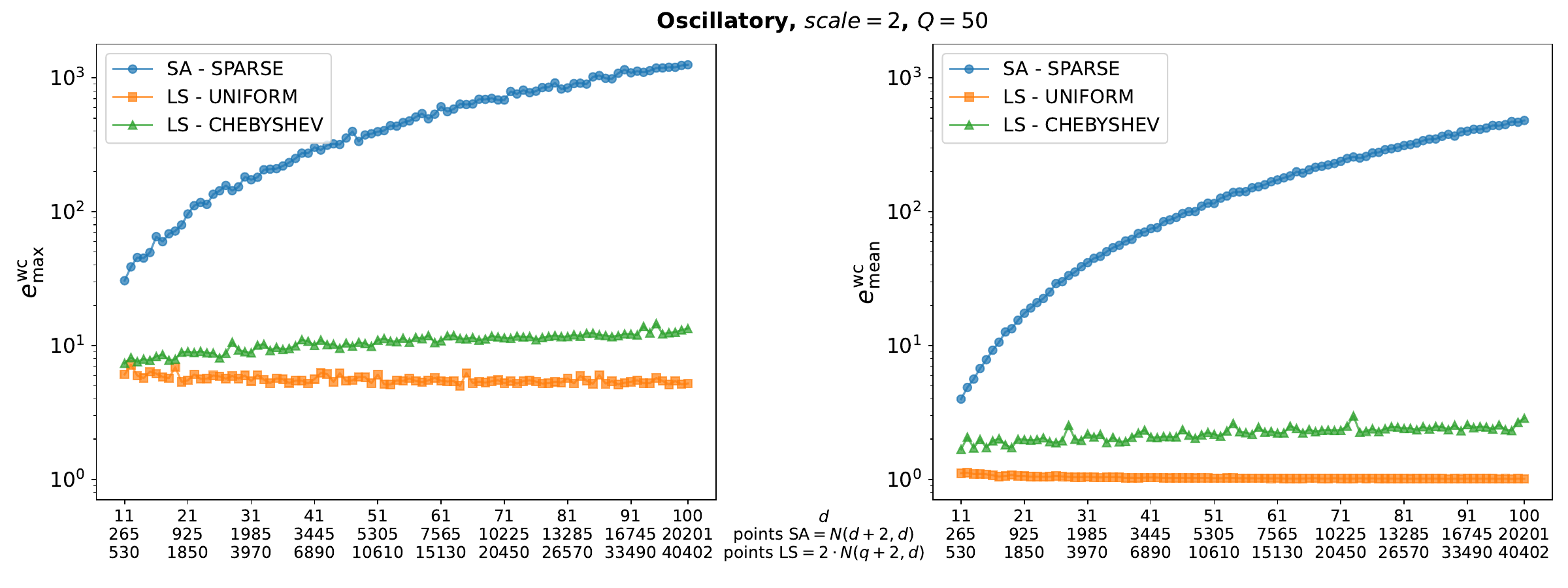}
    
    \caption{Some functions for fixed scale 1 and 2.}
	\label{fig:scale12_highdim_some_f_2}
\end{figure}

\begin{figure}[H]
	\centering
    \includegraphics[width=\widthoffigures\linewidth]{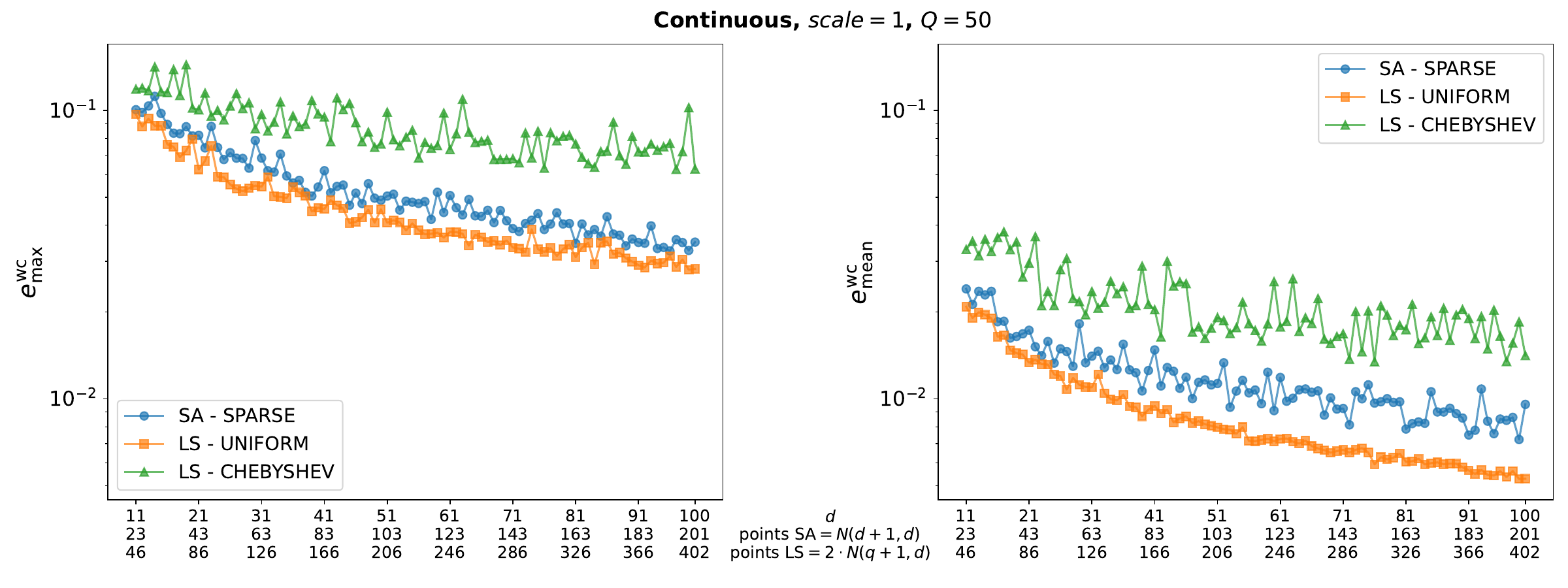}
    \includegraphics[width=\widthoffigures\linewidth]{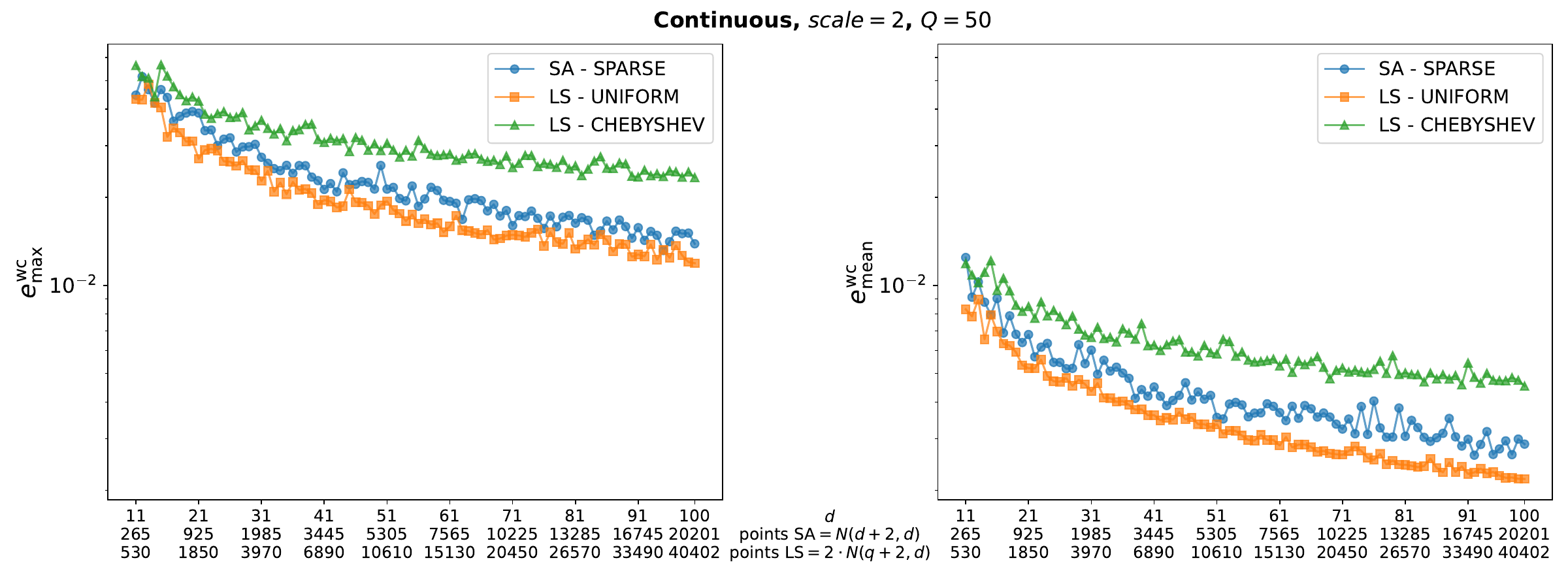}

    \includegraphics[width=\widthoffigures\linewidth]{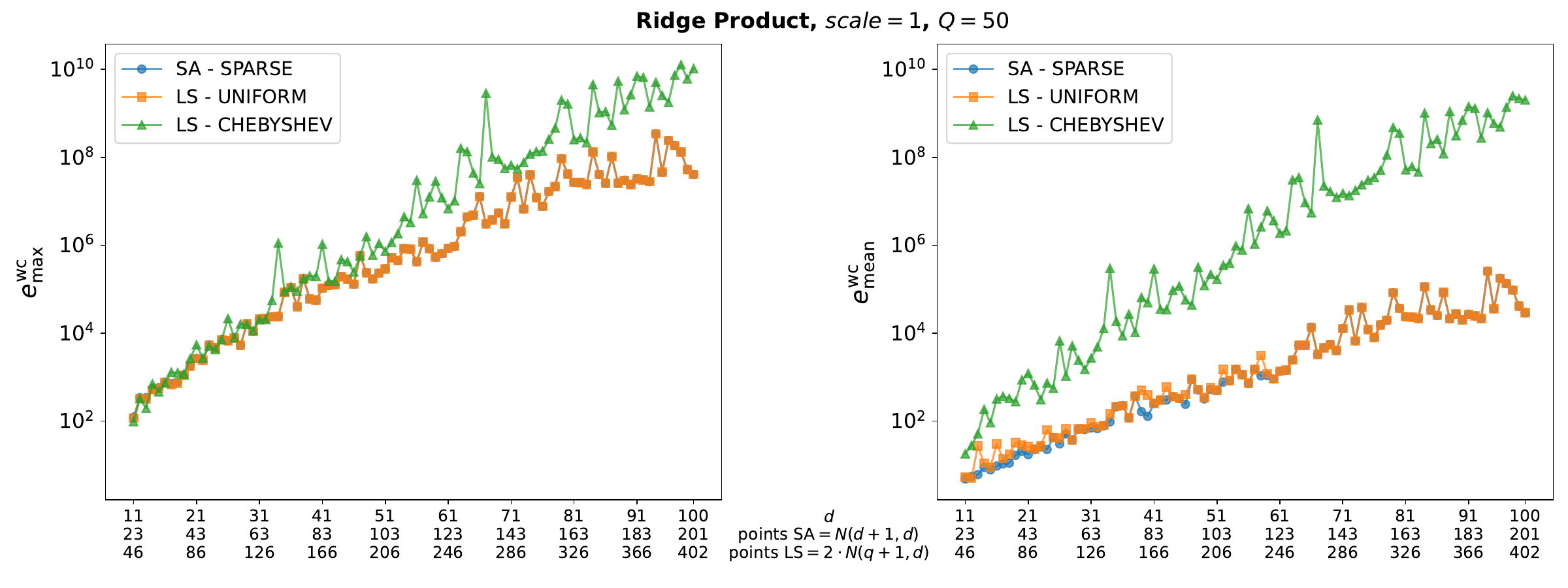}
    \includegraphics[width=\widthoffigures\linewidth]{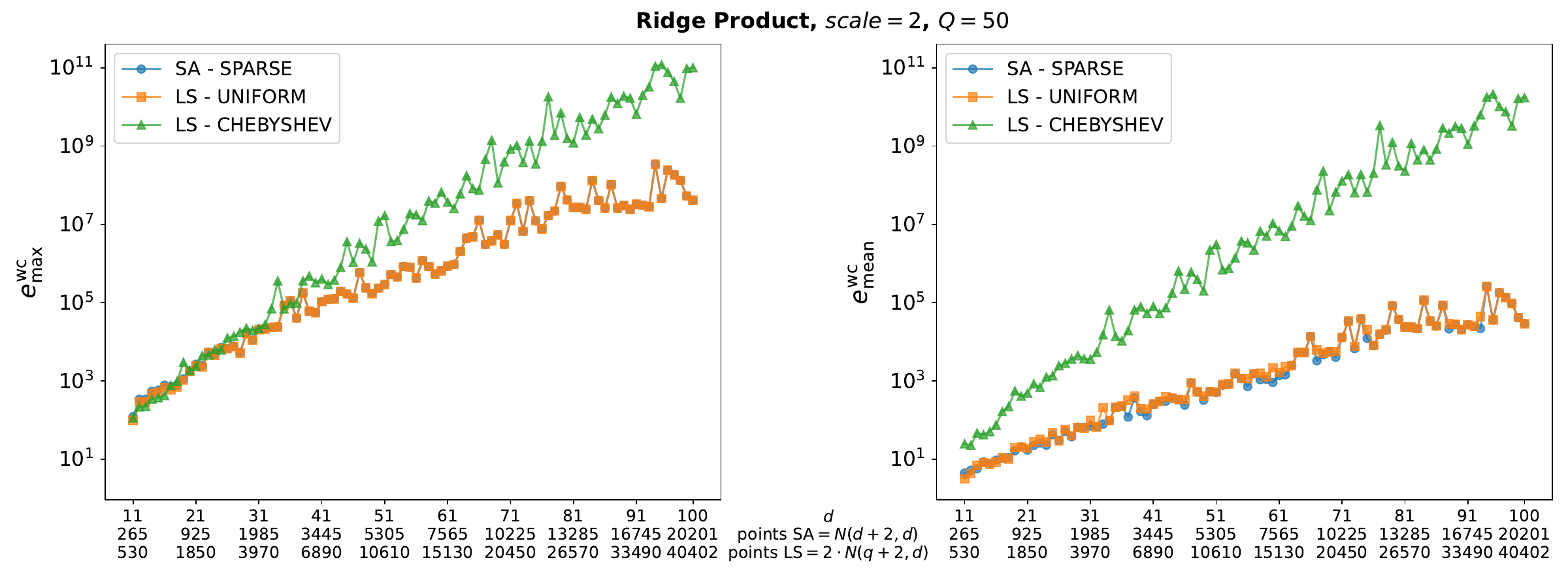}
    
    \caption{Some functions for fixed scale 1 and 2.}
	\label{fig:scale12_highdim_some_f_3}
\end{figure}

We now turn to our result in high dimensions up to 100. 
We can only consider small scales due to the fast increasing number of points. 
(See Figure~\ref{fig:runtime_comparison_high_dim} for runtimes.)
To have a decent amount of test-points, we use $M=100\cdot N(d+2, d)$ here.
Figures~\ref{fig:scale12_highdim_some_f_1} to \ref{fig:scale12_highdim_some_f_3} 
show that for many functions, 
Smolyak does not show a good performance compared to LS.
Moreover, LS-uniform often appears to be better than LS-Chebyshev.

Note that the large values of the error and the (non-)convergence indicated in the figures may be caused by a ``wrong'' scaling (in $d$) of the test functions. 
The result should therefore be understood, again, just as a comparison of the three methods. 


\medskip

\bigskip
\goodbreak

\subsection{Approximation of noise}
\label{sec:num-noise}

Finally, we made a few experiments related to \emph{stability} of the methods, 
i.e., how the methods behave under noisy data. 
%
For this, we test the algorithms 
against the $0$-function, 
but instead of actual function values, we provide 
Gaussian noise with mean $0$ and standard deviation $10^{-7}$. 

\begin{figure}[H]
	\centering
	\includegraphics[width=\widthoffigures\linewidth]{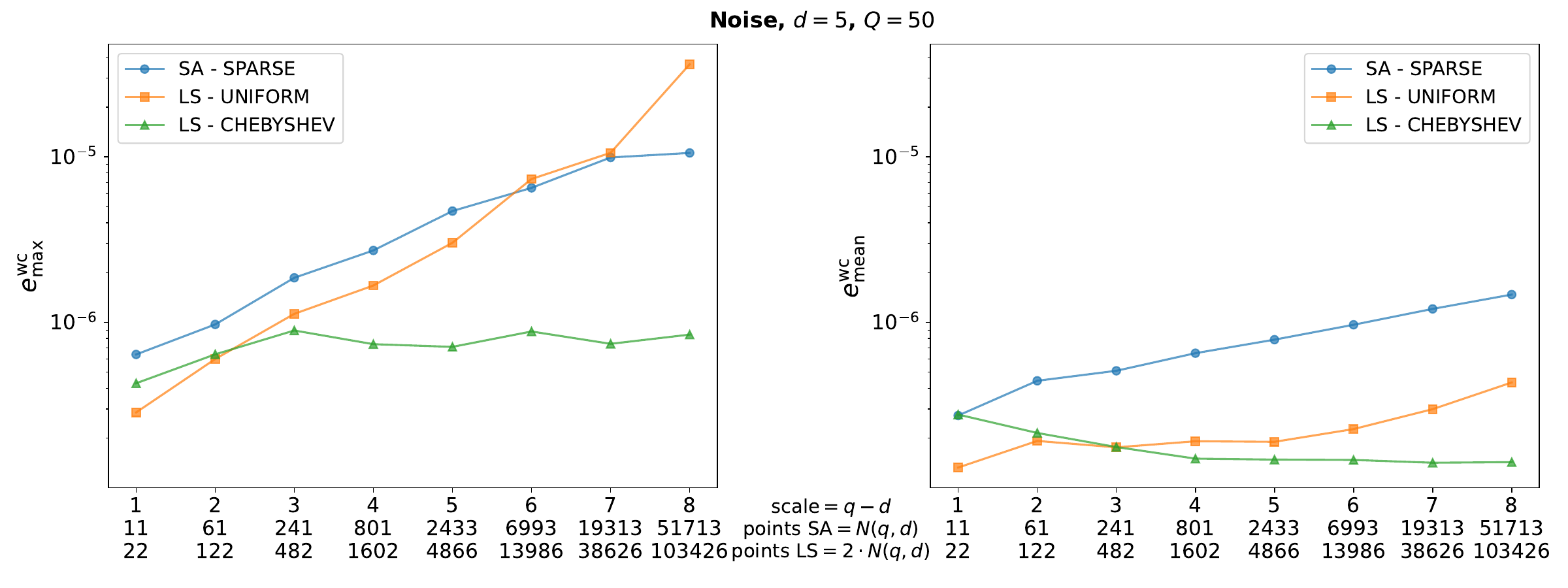}
    \includegraphics[width=\widthoffigures\linewidth]{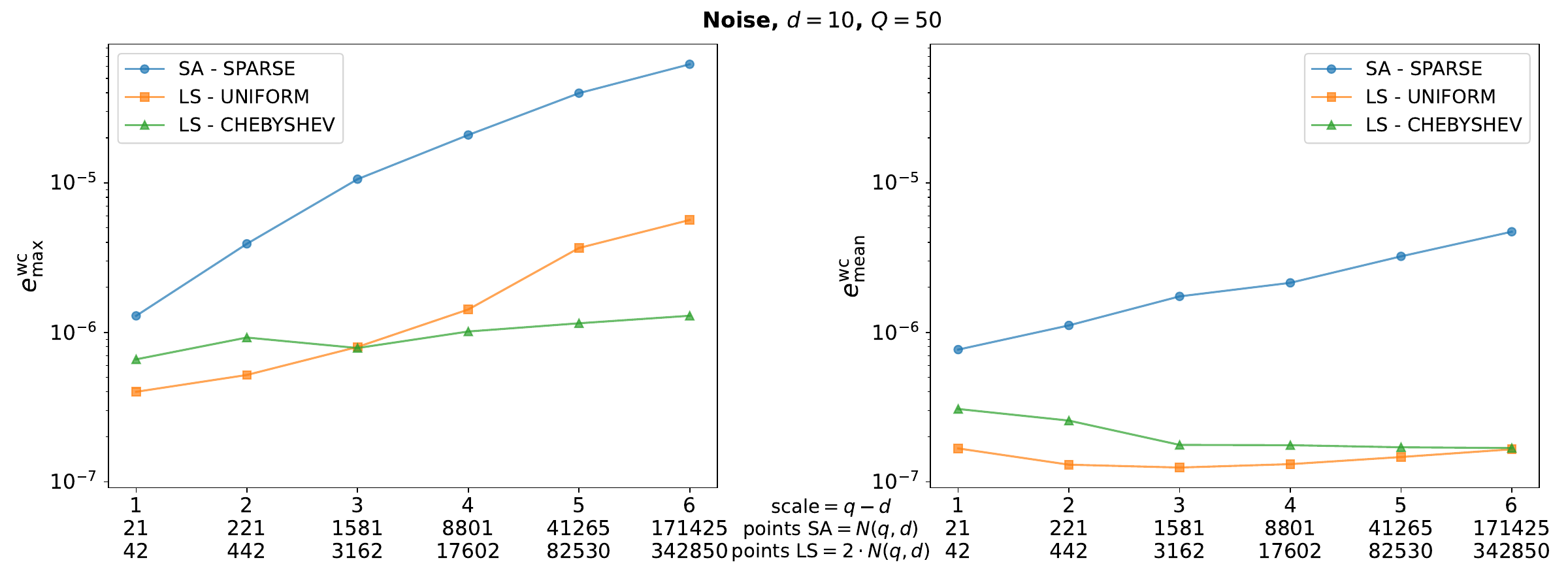}
	\caption{Noise for $d=5$ and $d=10$}
	\label{fig:noise_scale_510}
\end{figure}

It is known that interpolation is usually not a good idea in the case of noisy data.
One might expect that LS-Chebyshev is superior, at least asymptotically, 
and this can indeed be seen in Figure~\ref{fig:noise_scale_510}. 
We refer again to~\cite{CM17,MNST14} for a more detailed treatment.

In high dimensions, however,
LS-uniform seems to be more \emph{stable} than LS-Chebyshev.
This is to be understood theoretically. 

\begin{figure}[H]
	\centering
    \includegraphics[width=\widthoffigures\linewidth]{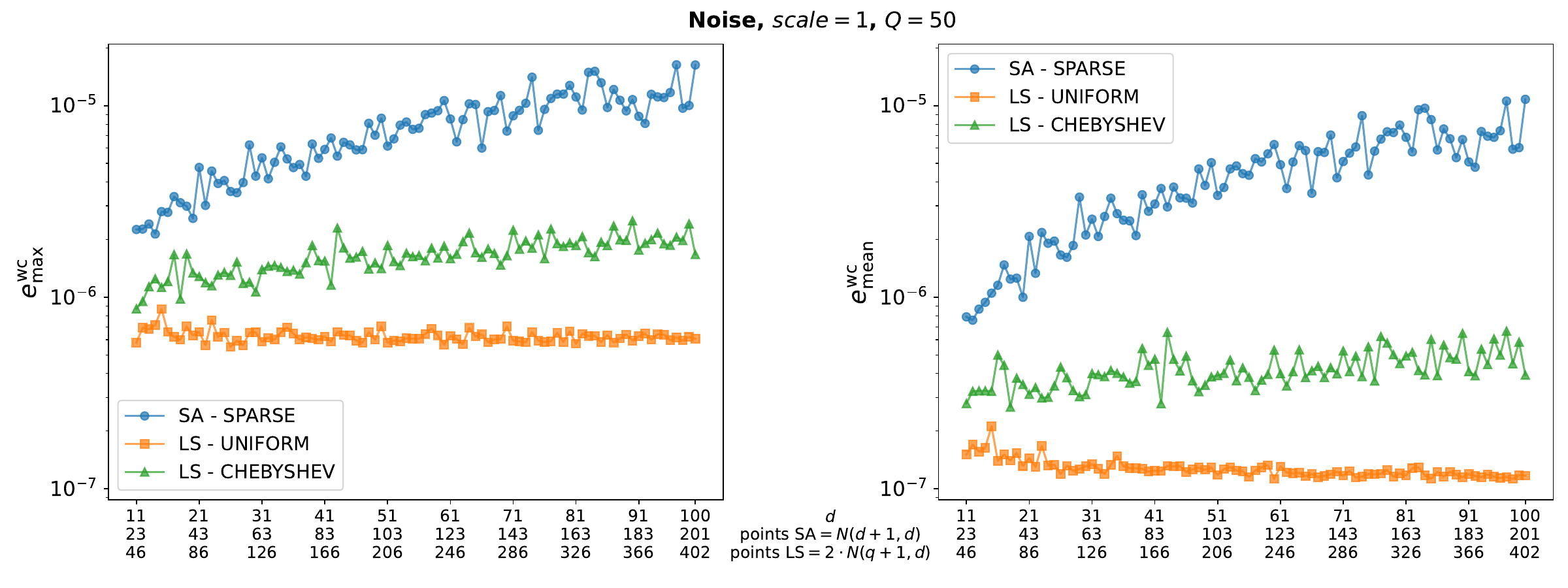}
	\includegraphics[width=\widthoffigures\linewidth]{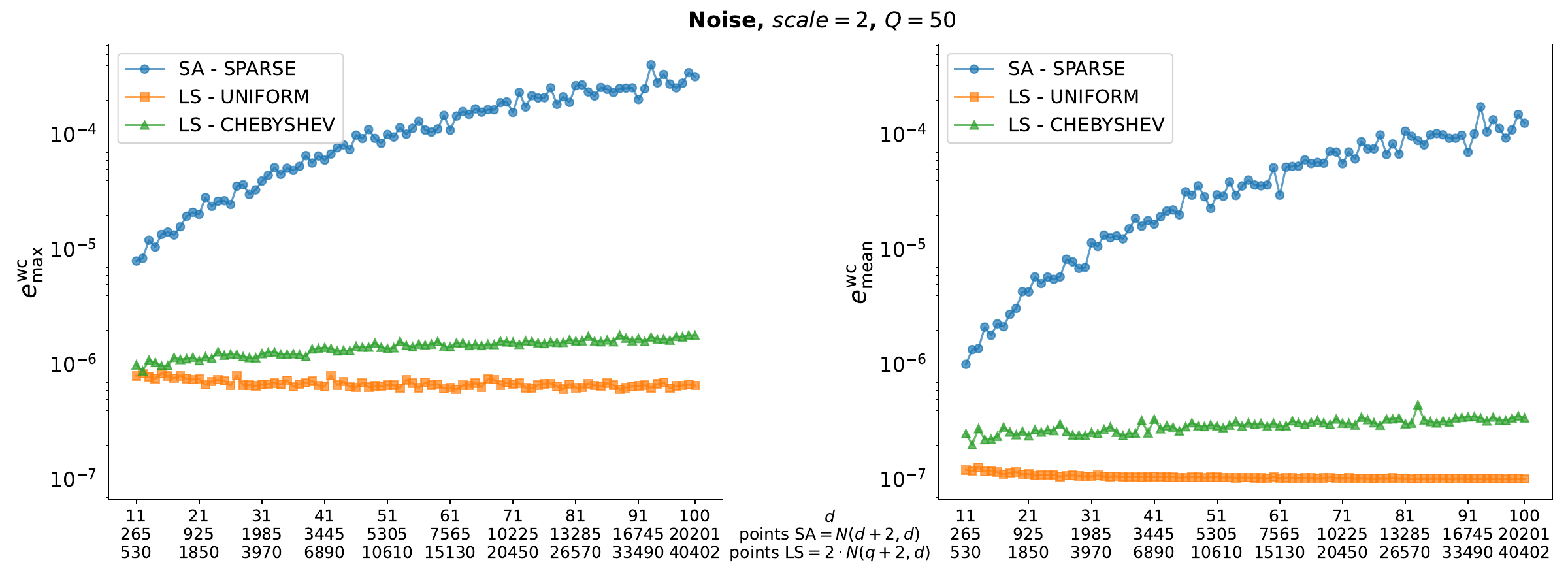}
	\caption{Noise for scales 1 and 2 in high dimensions.}
	\label{fig:noise_scale_12_high_dim}
\end{figure}


\section{Discussion}\label{sec:discussion}

Let us finally 
discuss some possible limitations 
and provide suggestions for future research.

The most important drawback of least squares is its running time and memory requirement: 
 Smolyak's algorithm is much faster than our implementation of LS. 
In particular, it was not possible to execute the latter, e.g., in dimension $d=10$ with $q=17$, i.e., scale 7, as this would have required 3TB of memory. For SA this is no problem at all.

\begin{figure}[H]
    \includegraphics[width=.95\linewidth]{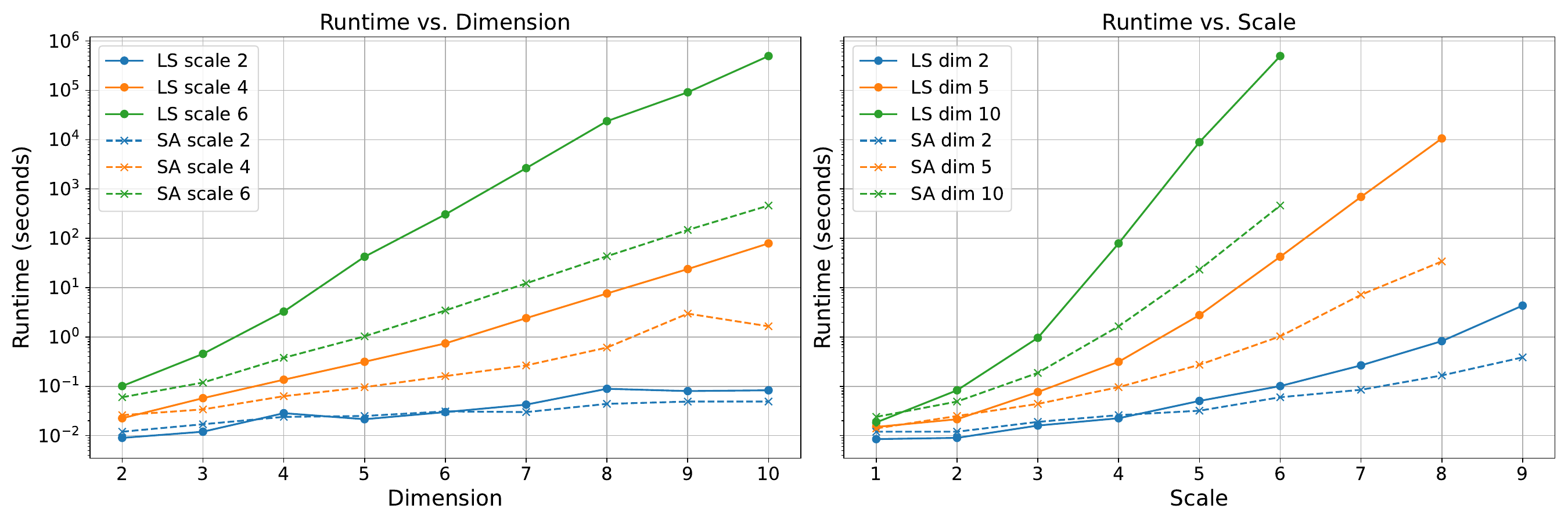}
    \captionof{figure}{Runtime for small dimensions}
    \label{fig:runtime_comparison_low_dim}
\end{figure}

\begin{figure}[H]
    \includegraphics[width=0.7\linewidth]{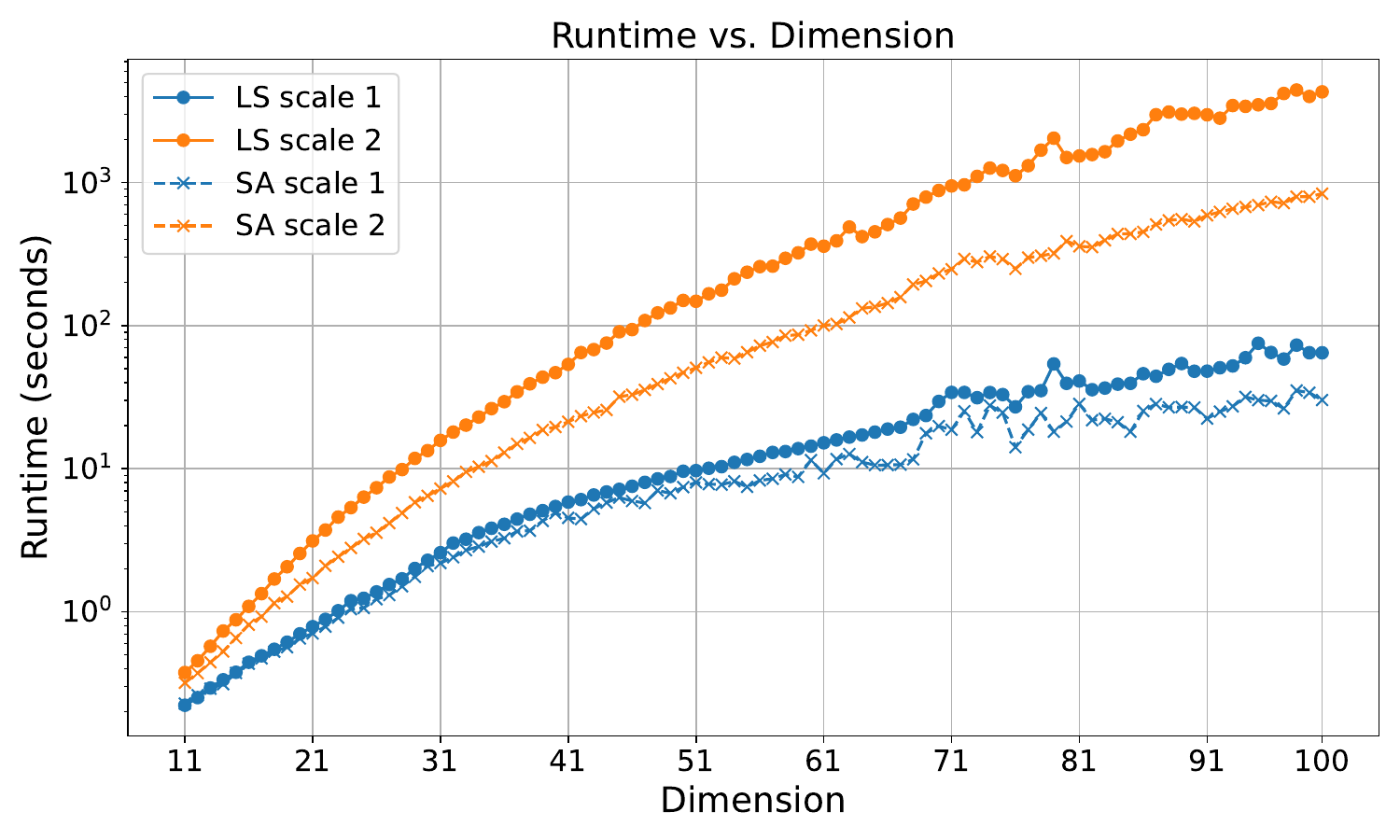}
    \captionof{figure}{Runtime for high dimensions}
    \label{fig:runtime_comparison_high_dim}
\end{figure}

Figures~\ref{fig:runtime_comparison_low_dim} and \ref{fig:runtime_comparison_high_dim} represent the runtime needed for our experiments for all algorithms and a selection of scales/dimensions. 
All of these were done with $Q=50$ realizations of each function class, and therefore with $500$ functions in parallel. The runtime of LS is the average runtime of LS-uniform and LS-Chebyshev.

\medskip

Here are some further comments:
\begin{itemize}
\itemsep=2mm
\item We used a rather straightforward implementation of the least squares method. 
We leave it for future research to perform numerical tests in the above setting with more elaborate methods to find an (approximate) solution to the huge linear system, like \emph{gradient descent} or \emph{batch learning}.

\item Motivated by~\cite{BNR00}, we used Smolyak's algorithm with a polynomial basis. Similarly, one may use piecewise linear basis functions, see e.g.~\cite{BG04} for details. 
Some small experiments did not show a difference in the reported results.
We did not try splines. 

\item Similarly, we might consider different \emph{approximation spaces} $V$ for least squares. 
We used the same space as Smolyak's algorithm for a \emph{fair} comparison. It would be of interest to find and test other natural choices.

\item We also tried to benchmark the three algorithms by using test-points sampled from the Chebyshev density. This led to unessential differences in the estimated errors which 
we do not discuss here, because it corresponds to a different error criterion. The results can be found on GitHub.

\item Visualizations of the errors using boxplots (indicating other statistical effects) can be found on GitHub.
We did not include them in this paper, because they were visually unappealing and did not show anything interesting (for us). 

\end{itemize}

\bigskip
\noindent \textbf{Acknowledgement.}
%
We thank Erich Novak for comments that led to some of our numerical tests and to improvements in presentation.
We also thank Manuel Kauers from the Institute of Algebra at JKU for letting us perform the computations reported in this paper on the computers of his group. 

JE and EM are supported by the FWF Start Project ``Quantum shadows: scalable quantum-to-classical converters'' (Richard K\"ung) 
and by JKU as student research employees.
MU is supported by the Austrian Federal Ministry of Education, Science and Research via the Austrian Research Promotion Agency (FFG) through the project FO999921407 (HDcode) funded by the European Union via NextGenerationEU.

\medskip

\linespread{0.9}







\begin{thebibliography}{99}

\bibitem{ABB99}
E. Anderson, Z. Bai, C. Bischof et al.,
{\emph{{LAPACK} users' guide}},
Society for Industrial and Applied Mathematics, 3, 1999,
\doi{10.1137/1.9780898719604}.


\bibitem{BSU23}
F. Bartel, M. Sch\"afer and T. Ullrich, 
{\emph{Constructive subsampling of finite frames
with applications in optimal function recovery}}, Applied Comput. Harmon.
Anal., 65, 209–-248, 2023,
\doi{10.1016/j.acha.2023.02.004}.


\bibitem{BNR00}
V. Barthelmann, E. Novak and K. Ritter,
{\emph{High dimensional polynomial interpolation on sparse grids}},
Advances in Computational Mathematics 12,  273--288, 2000, 
\doi{10.1023/A:1018977404843}.


\bibitem{BCLSV11}
L. Bos, J.P. Calvi, N. Levenberg, A. Sommariva and M. Vianello,
{\emph{Geometric weakly admissible meshes, discrete least squares 
approximations and approximate Fekete points}},
Math. Comp., 80, 1623--1638, 2011,
\doi{10.1090/S0025-5718-2011-02442-7}.


\bibitem{Bos18a}
L. Bos,
{\emph{Fekete points as norming sets}},
Dolomites Research Notes on Approximation, 11(4), 26-–34, 2018,
\doi{10.14658/PUPJ-DRNA-2018-4-3}.


\bibitem{Bos18b}
L. Bos,
{\emph{On optimal designs for a d-cube}},
Dolomites Research Notes on Approximation, 15(4), 20--34, 2018,
\doi{10.14658/PUPJ-DRNA-2022-4-3}.


\bibitem{BPV19}
L. Bos, F. Piazzon and M. Vianello,
{\emph{Near optimal polynomial regression on norming meshes}},
Sampling Theory and Applications, 2019,
\doi{10.1109/SampTA45681.2019.9030910}.




\bibitem{BG04}
H.J. Bungartz and M. Griebel,
{\emph{Sparse grids}},
Acta Numerica 13, 147--269, 2004, 
\href{https://www.cambridge.org/core/journals/acta-numerica/article/sparse-grids/47EA2993DB84C9D231BB96ECB26F615C}{doi: 10.1017/S0962492904000182}


\bibitem{BKSVZ24}
L. Becker, O. Klein, J. Slote, A. Volberg and H. Zhang,
{\emph{Dimension-free discretizations of the uniform norm by small product sets}},
Invent. math., 239, 469--503, 2025,
\doi{10.1007/s00222-024-01306-9}.






\bibitem{CL08}
J.P. Calvi and N. Levenberg,
{\emph{Uniform approximation by discrete least squares polynomials}},
J. Approx. Theory, 152(1), 82--100, 2009, 
\doi{10.1016/j.jat.2007.05.005}.


\bibitem{CCMNT15}
A. Chkifa, A. Cohen, G. Migliorati, F. Nobile and R. Tempone,
{\emph{Discrete least squares polynomial approximation with random evaluations 
- application to parametric and stochastic elliptic PDEs}},
M2AN, 49, 815–-837, 2015,
\doi{10.1051/m2an/2014050}.


\bibitem{CD23}
A. Chkifa and M. Dolbeault, 
{\it Randomized least-squares with minimal oversampling and
interpolation in general spaces}, 
SIAM J. Numer. Anal.  62(4), 1515--1538, 2024,
\doi{10.1137/23M160178X}.


\bibitem{CDL13}
A. Cohen, M.A. Davenport and D. Leviatan,
{\emph{On the stability and accuracy of least squares approximations}},
Foundations of Computational Mathematics, 13, 819--834, 2013,
\doi{10.1007/s10208-013-9142-3}.


\bibitem{CM17}
A. Cohen, G. Migliorati,
{\emph{Optimal weighted least squares methods}},
SMAI J. Comput. Math. 3, 181–-203, 2017,
\doi{10.5802/smai-jcm.24}.


\bibitem{CS13}
C. Coleman and S. Lyon,
{\emph{Efficient implementations of Smolyak's algorithm for function approximation in Python and Julia}}
2013, \url{https://github.com/EconForge/Smolyak}.




\bibitem{DC22}
M. Dolbeault and A. Cohen,
{\emph{Optimal pointwise sampling for $L_2$ approximation}}, 
J. Complexity, 68, 101602, 2022,
\doi{10.1016/j.jco.2021.101602}.


\bibitem{DKU} M.~Dolbeault, D.~Krieg and M.~Ullrich,  
{\it A sharp upper bound for sampling numbers in $L_2$}, 
{Appl. Comput. Harmon. Anal.} 63, 113--134, 2023,
\doi{10.1016/j.acha.2022.12.001}.




\bibitem{DP23}
F. Dai and A. Prymak,
{\emph{Optimal polynomial meshes exist on any multivariate convex domain}},
Found Comput Math, 24, 989--1018, 2024,
\doi{10.1007/s10208-023-09606-x}.


\bibitem{DTU18}
D. D{\~u}ng, V. Temlyakov and T. Ullrich,
{\emph{Hyperbolic cross approximation}}
Birkh\"auser Cham, 2018,
\doi{10.1007/978-3-319-92240-9}.


\bibitem{EZ66}
H. Ehlich and K. Zeller,
{\emph{Auswertung der Normen von Interpolationsoperatoren}},
Mathematische Annalen, 105--112, 1966,
\doi{10.1007/BF01429047}.


\bibitem{Gen84}
A. Genz, 
{\emph{Testing multidimensional integration routines}},
Proc. of International Conference on Tools, Methods and Languages for 
Scientific and Engineering Computation, 81--94, 1984.


\bibitem{Gen87}
A. Genz, 
{\emph{A package for testing multiple integration subroutines}},
Numerical integration: Recent developments, software and applications, 203,
337--340, 1987,
\doi{10.1007/978-94-009-3889-2\_33}.

\bibitem{GHM25}
M.~Griebel, H.~Harbrecht and M.~Multerer, 
\emph{Kernel interpolation on sparse grids}, 
arXiv:2505.12282, 2025.

\bibitem{Gr19}
K. Gröchenig,
{\it Sampling, Marcinkiewicz-Zygmund inequalities, approximation, and quadrature rules}
, J. Approx. Theory, 257:105455, 2020,
\doi{10.1016/j.jat.2020.105455}.


\bibitem{GNZ20}
L. Guo, A. Narayan and T. Zhou,
{\emph{Constructing least-squares polynomial approximations}},
SIAM Review, 62(2), 483--508, 2020,
\doi{10.1137/18M1234151}.


\bibitem{HMW20}
C.R. Harris, K.J.Millman and S. van der Walt et al.,
{\emph{Array programming with {NumPy}}},
Nature, 585, 357--362, 2020,
\doi{10.1038/s41586-020-2649-2}.


\bibitem{Jo92}
M.L. Johnson, 
{\emph{Why, when, and how biochemists should use least squares},
Analytical biochemistry, 206(2), 215-225, 1992,
\doi{10.1016/0003-2697(92)90356-c}


\bibitem{JMMV14}
K.L. Judd, L. Maliar, S. Maliar and R. Valero,
{\emph{Smolyak method for solving dynamic economic models: Lagrange 
interpolation, anisotropic grid and adaptive domain}}
J. Economic Dynamics and Control, 44, 92--123, 2014,
\doi{10.1016/j.jedc.2014.03.003}.


\bibitem{KKT23}
B.~Kashin, S.~Konyagin and V.~Temlyakov, 
{\emph{Sampling discretization of the uniform norm}},
Constr. Approx., 57(2), 663--694, 2023,
\doi{10.1007/s00365-023-09618-4}


\bibitem{KKLT} B. Kashin, E. Kosov, I. Limonova and V. Temlyakov,
{\it Sampling discretization and related problems},
J. Complexity, 101653, 2022.


\bibitem{Kroo11}
A. Kro{\'o},
{\emph{On optimal polynomial meshes}},
J. Approx Theory, 163(9), 1107--1124, 2011,
\doi{10.1016/j.jat.2011.03.007}.



\bibitem{KPUU25-uniform}
D.~Krieg, K.~Pozharska, M.~Ullrich and T.~Ullrich,
{\emph{Sampling projections in the uniform norm}}
to appear in Forum Math Sigma, 2025,
\url{https://arxiv.org/abs/2401.02220}.

\bibitem{KPUU25-general}
D.\,Krieg, K.\,Pozharska, M.\,Ullrich, T.\,Ullrich,
{\it Sampling recovery in $L_2$ and other norms}, 
to appear in Math. Comp., 2023,
\url{https://arxiv.org/abs/2305.07539}.

\bibitem{KSUW} D. Krieg, P. Siedlecki, M. Ullrich and H. Wo\'zniakowski,
{\it Exponential tractability of $L_2$-approximation with function values},
Adv. Comput. Math. 49, 18, 2023,
\doi{10.1007/s10444-023-10021-7}

\bibitem{KU21}
D.~Krieg and M.~Ullrich,
\emph{Function values are enough for \({L}_2\)-approximation},
Found. Comp. Mathematics, 21(4), 1141--1151, 2021,
\doi{10.1007/s10208-020-09481-w}.

\bibitem{KU21b} D.~Krieg and M.~Ullrich,
{\it Function values are enough for $L_2$-approximation: Part II}, 
J. Complexity 66, 101569, 2021,
\doi{10.1016/j.jco.2021.101569}.



\bibitem{LT84}
L.R. Lines and S. Treitel,
{\emph{A review of least‐squares inversion and its application to geophysical problems}},
Geophysical prospecting, 32(2), 159-186, 1984. 

\bibitem{MNST14}
G. Migliorati, F. Nobile, E. von Schwerin and R. Tempone,
{\emph{Analysis of discrete projection on polynomial spaces with random evaluations}},
Found Comput Math, 14, 419-–456, 2014,
\doi{10.1007/s10208-013-9186-4}.



\bibitem{NJZ17}
A.~Narayan, J.~Jakeman, T.~Zhou, 
\emph{A Christoffel function weighted least squares algorithm for collocation approximations},
Math. Comp. 86 (306), 1913--1947, 2017,
\url{http://arxiv.org/abs/1412.4305}.



\bibitem{NR96}
E. Novak and K. Ritter,
{\emph{High dimensional integration of smooth functions over cubes}},
Num. Mathematik, 75, 79--97, 1996,
\doi{10.1007/s002110050231}.


\bibitem{NR99}
E. Novak and K. Ritter,
{\emph{Simple cubature formulas with high polynomial exactness}},
Constr. Approx., 15, 499–-522, 1999,
\doi{10.1007/s003659900119}.


\bibitem{NX23}
A. Narayan and Y. Xu,
{\emph{Randomized weakly admissible meshes}},
J. Approx. Theory, 285, 105835, 2023,
\doi{10.1016/j.jat.2022.105835}.


\bibitem{Reichel86}
L. Reichel
{\emph{On polynomial approximation in the uniform norm by the discrete least squares method}},
BIT, 26(3), 349--368, 1986,
\doi{10.1007/BF01933715}.


\bibitem{Smo63}
S.A.~Smolyak,
{\emph{Quadrature and interpolation formulas for tensor products of certain classes of functions}},
Dokl. Akad. Nauk SSSR, 148(5), 1042--1045, 1963,
\url{http://mathscinet.ams.org/mathscinet-getitem?mr=0147825}


\bibitem{SU23}
M. Sonnleitner, M. Ullrich,
{\emph{On the power of iid information for linear approximation}},
J. Appl. Numer. Anal., 1, 88--126, 2023,
\doi{10.30970/ana.2023.1.88}.

\bibitem{Sto15}
M.K. Stoyanov,
{\emph{User manual: {TASMANIAN} sparse grids}},
ORNL/TM-2015/596, 2015.


\bibitem{SLBM13}
M.K. Stoyanov, D. Lebrun-Grandie, J. Burkardt and D. Munster,
{\emph{{Tasmanian}}},
2013,
\doi{10.11578/dc.20171025.on.1087}.






\bibitem{SB13}
S. Surjanovic and D. Bingham,
{\emph{Virtual library of simulation experiments: Test functions and datasets}},
2013, retrieved March 29, 2025, \url{https://www.sfu.ca/~ssurjano/integration.html}.


\bibitem{T21}
V.\,N.\,Temlyakov,
{\emph{On optimal recovery in $L_2$}},  
J. Complexity, 65, 101545, 2021,
\doi{10.1016/j.jco.2020.101545}.


\bibitem{Tem86}
V.N. Temlyakov,
{\emph{Approximation of periodic functions of several variables by 
trigonometric polynomials, and widths of some classes of functions}},
Math. USSR-Izv., 27(2), 285--322, 1986,
\doi{10.1070/IM1986v027n02ABEH001179}.


\bibitem{Ull20}
M. Ullrich,
{\emph{On the worst-case error of least squares algorithms for 
$L_2$-approximation with high probability}},
J. Complexity, 101484, 2020,
\doi{10.1016/j.jco.2020.101484}.




\bibitem{VGO20}
P. Virtanen, R. Gommers and T.E. Oliphant et al.,
{\emph{{SciPy} 1.0: Fundamental algorithms for scientific computing in Python}},
Nature Methods, 17, 261--272, 2020,
\doi{10.1038/s41592-019-0686-2}.


\bibitem{WW95}
G.W. Wasilkowski and H. Wo\'zniakowski,
{\emph{Explicit cost bounds of algorithms for multivariate tensor product problems}},
J. Complexity, 11(1), 1--56, 1995,
\doi{10.1006/jcom.1995.1001}.





}
\end{thebibliography}
\end{document}